\documentclass[12pt,reqno]{amsart}

\usepackage{amsfonts,color,amsthm,amsmath,amssymb}

\usepackage{color}
\allowdisplaybreaks[3]

\def\Box{\vcenter{\vbox{\hrule\hbox{\vrule
     \vbox to 8.8pt{\hbox to 10pt{}\vfill}\vrule}\hrule}}}

\def\qed{{\hfill$\square$}}
\def\proof{{\vspace{-0.0cm}\bf Proof: \,}}

\def\Z{{\mathbb Z}}

\def\F{{\mathbb F}}

\def\mod{{\mathrm{mod\,\,}}}

\def\Tr{{\mathrm{Tr}}}

\def\PG{{\mathrm{PG}}}

\newtheorem{theorem}{Theorem}[section]

\newtheorem{lemma}[theorem]{Lemma}
\newtheorem{remark}[theorem]{Remark}
\newtheorem{corollary}[theorem]{Corollary}

\newtheorem{proposition}[theorem]{Proposition}

\numberwithin{equation}{section}

\begin{document}
\title[Menon-Hadamard difference sets]{Generalized constructions of Menon-Hadamard difference sets
}

\author[Momihara and Xiang]{Koji Momihara$^{\dagger}$ and Qing Xiang$^{\ast}$}
\thanks{$^{\dagger}$
Koji Momihara was supported by 
JSPS under Grant-in-Aid for Young Scientists (B) 17K14236 and Scientific Research (B) 15H03636.}

\thanks{$^{\ast}$
Qing Xiang was supported by an NSF grant DMS-1600850.}

\address{Faculty of Education, Kumamoto University, 2-40-1 Kurokami, Kumamoto 860-8555, Japan} \email{momihara@educ.kumamoto-u.ac.jp}

\address{Department of Mathematical Sciences, University of Delaware, Newark, DE 19716, USA} \email{xiang@math.udel.edu}

\begin{abstract} We revisit the problem of constructing Menon-Hadamard difference sets. In 1997, Wilson and Xiang gave a general framework for constructing Menon-Hadamard difference sets by using a combination of a spread and four projective sets of type Q in ${\mathrm{PG}}(3,q)$. They also found examples of suitable spreads and projective sets of type Q for $q=5,13,17$.  Subsequently, Chen (1997) succeeded in finding a spread and four projective sets of type Q in ${\mathrm{PG}}(3,q)$ satisfying the conditions in the Wilson-Xiang construction for all odd prime powers $q$. Thus, he showed that there exists a Menon-Hadamard difference set of order $4q^4$ for all odd prime powers $q$. However, the projective sets of type Q found by Chen have automorphisms different from those of the examples constructed by Wilson and Xiang.  In this paper, we first generalize Chen's construction of projective sets of type Q by using ``semi-primitive'' cyclotomic classes. This demonstrates that the construction of projective sets of type Q satisfying the conditions in the Wilson-Xiang construction is much more flexible than originally thought. Secondly, we give a new construction of spreads and projective sets of type Q in ${\mathrm{PG}}(3,q)$ for all odd prime powers $q$, which generalizes the examples found by Wilson and Xiang. This solves a problem left open in Section 5 of the Wilson-Xiang paper from 1997.
\end{abstract}

\maketitle
\section{Introduction}

Let $G$ be an additively written abelian group of order $v$. A $k$-subset $D$ of $G$ is called 
a {\it $(v,k,\lambda)$ difference set} if the list of differences ``$x-y$,  $x,y\in D, x\neq y$", represents each nonidentity element of $G$ exactly $\lambda$ times. 
In this paper, we revisit the problem of constructing Menon-Hadamard difference sets, namely those difference sets with parameters $(v,k,\lambda)=(4m^2,2m^2-m,m^2-m)$, where $m$ is a positive integer. It is well known that a Menon-Hadamard difference set generates a regular Hadamard matrix of order $4m^2$. So by contructing Menon-Hadamard difference sets in groups of order $4m^2$, we obtain regular Hadamard matrices of order $4m^2$.

The main problem in the study of Menon-Hadamard difference sets is: For each positive integer $m$, which groups of order $4m^2$ contain a Menon-Hadamard difference set. We give a brief survey of results on this problem in the case where the group under consideration is abelian. First we mention a product theorem of Turyn~\cite{T84}: If there are Menon-Hadamard difference sets in abelian groups $H\times G_1$ and $H\times G_2$, respectively, where $|H|=4$ and $|G_i|$, $i=1,2$, are squares, then there also exists a Menon-Hadamard difference set in $H\times G_1\times G_2$. With Turyn's product theorem in hand, in order to construct Menon-Hadamard difference sets, one should start with the case where the order of the abelian group is $4q$ with $q$ an even power of a prime. In the case where $q$ is an even power of $2$, that is, $G$ is an abelian $2$-group, the existence problem was completely solved in \cite{K93} after much work was done in \cite{davis}; it was shown that there exists a Menon-Hadamard difference set in an abelian group $G$ of order $2^{2t+2}$ if and only if the exponent of $G$ is less than or equal to $2^{t+2}$.

In the case where $q$ is an even power of an odd prime, Turyn \cite{T84} observed that there exists a Menon-Hadamard difference set in $H\times (\Z_3)^2$; hence by the product theorem, there is a Menon-Hadamard difference set in $H\times (\Z_3)^{2t}$ for any positive integer $t$. On the other hand, McFarland \cite{mcfarland} proved that if an abelian group of order $4p^2$, where $p$ is a prime, contains a Menon-Hadamard difference set, then $p=2$ or $3$. After McFarland's paper \cite{mcfarland} was published, it was conjectured \cite[p.~287]{jung} that if an abelian group of order $4m^2$ contains a Menon-Hadamard difference set, then $m=2^r3^s$ for some nonnegative integers $r$ and $s$.  So it was a great surprise when Xia~\cite{X} constructed a Menon-Hadamard difference set in $H\times \Z_{p}^4$ for any odd prime $p$ congruent to $3$ modulo $4$. Xia's method of contruction depends on very complicated computations involving cyclotomic classes of finite fields; it  was later simplified by Xiang and Chen~\cite{XC} by using a character theoretic approach. Moreover, in \cite{XC}, the authors also asked whether a certain family of 3-weight projective linear code exists or not, since such projective linear codes are needed for the construction of Menon-Hadamard difference set in the group $H\times (\Z_p)^4$, where $p$ is a prime congruent to 1 modulo 4.

Van Eupen and Tonchev~\cite{ET} found the required 3-weight projective linear codes when $p=5$, hence constructed Menon-Hadamard difference sets in $\Z_2^2\times \Z_5^4$, which are the first examples of abelian Menon-Hadamard difference sets in groups of order $4p^4$, where $p$ is a prime congruent to $1$ modulo $4$. Inspired by these examples, Wilson and Xiang~\cite{WX97} gave a general framework for constructing Menon-Hadamard difference sets in the groups $H\times G$, where $H$ is either group of order $4$ and $G$ is an elementary abelian group of order $q^4$,  $q$ an odd prime power, using a combination of a spread and four projective sets of type Q in $\PG(3,q)$. (See Section~\ref{sec:twonew} for the definition of projective sets of type Q.)
Wilson and Xiang \cite{WX97} also found examples of suitable spreads and the required projective sets of type Q when $q=5,13,17$. They used  $\F_{q^2}\times \F_{q^2}$ as a model of the four-dimensional vector space $V(4,q)$ over $\F_q$, and considered projective sets of type Q with the automorphism 
\[
T'=\begin{pmatrix}
\omega^2  &0 \\
0 &\omega^{-2}
\end{pmatrix}, 
\]
where $\omega$ is a primitive element of $\F_{q^2}$. However, the existence of the required projective sets of type Q with this prescribed automorphism remained unsolved for  $q>17$.
 
Immediately after \cite{WX97} appeared, Chen~\cite{Ch97} succeeded in showing the existence of a combination of a spread and four projective sets of type Q in $\PG(3,q)$ satisfying the conditions in the Wilson-Xiang construction for all odd  prime powers $q$. As a consequence, Chen \cite{Ch97} obtained the following theorem by applying Turyn's product theorem in \cite{T84}. 
\begin{theorem}\label{thm:Chen}
Let $p_i$, $i=1,2,\ldots,s$, be odd primes and $t_i$, $i=1,2,\ldots,s$, be positive integers. Furthermore, let $H$ be either group of order $4$ and  $G_i$, $i=1,2,\ldots,s$, be an elementary abelian group of order $p_i^{4t_i}$. Then, 
there exists a Menon-Hadamard difference set in $H\times G_1\times G_2\times \cdots \times G_s$. 
\end{theorem}
Here, Chen \cite{Ch97} found projective sets of type Q in $\PG(3,q)$ with the following automorphism
\[
T=\begin{pmatrix}
\omega^2  &0 \\
0 &\omega^{2} 
\end{pmatrix},
\]
which is obviously  
different from that of the projective sets of type Q found by Wilson and Xiang \cite{WX97}.  Thus, the existence problem of projective sets of type Q in $\PG(3,q)$ with the prescribed automorphism $T'$ remained open. 

The objectives of this paper are two-fold. First, we give a generalization of Chen's construction of projective sets of type Q by using ``semi-primitive'' cyclotomic classes. This demonstrates that the construction of projective sets of type Q satisfying the conditions in the Wilson-Xiang construction is much more flexible than originally thought. In particular, the proof of the candidate sets are projective sets of type Q is much simpler than that in \cite{Ch97}. 
Second, we show the existence of a combination of a spread and four projective sets of type Q with automorphism $T'$ for all odd prime powers $q$. Our construction generalizes the examples found by Wilson and Xiang in \cite{WX97}; this solves the problem left open in Section 5 of \cite{WX97}. 



\section{Preliminaries}

\subsection{Characters of finite fields}
In this subsection, we collect some auxiliary results on characters of finite fields. We 
assume that the reader is familiar with basic theory of characters of finite fields as in \cite[Chapter~5]{LN97}. 

Let $p$ be a prime and $s,f$ be positive integers. We set $q=p^s$, and denote the finite field of order $q$ by $\F_{q}$. 
Let $\Tr_{q^f/q}$ be the trace map from $\F_{q^f}$ to $\F_{q}$, which is  defined by 
\[
\Tr_{q^f/q}(x)=x+x^q+\cdots+x^{q^{f-1}}, \quad x\in \F_{q^f}. 
\] 

Let $\omega$ be a fixed primitive element of $\F_q$, $\zeta_p$ a fixed (complex) primitive $p$th root of unity, and $\zeta_{q-1}$ a (complex) $q-1$th root of unity.  
The character $\psi_{\F_{q}}$ of the additive group of $\F_{q}$ defined by $\psi_{\F_{q}}(x)=\zeta_p^{\Tr_{q/p}(x)}$, $x\in \F_q$, is called the {\it canonical 
additive character}  of $\F_{q}$. Then, each additive character is given by $\psi_a(x)=\psi_{\F_q}(ax)$, $x\in \F_{q}$, where $a\in \F_q$.  
On the other hand, each multiplicative character is given by 
$\chi^j(\omega^\ell)=\zeta_{q-1}^{j\ell}$, 
$\ell=0,1,\ldots,q-2$, where 
$j=0,1,\ldots,q-2$. 

For a multiplicative character $\chi$ of $\F_q$, the character sum defined by 
\[
G_q(\chi)=\sum_{x\in \F_q^\ast}\chi(x)\psi_{\F_q}(x)
\]
is called a {\it Gauss sum} of $\F_q$. Gauss sums satisfy the following basic properties: (1) $G_q(\chi)\overline{G_q(\chi)}=q$ if $\chi$ is nontrivial; 
(2) $G_q(\chi^{-1})=\chi(-1)\overline{G_q(\chi)}$; (3) 
$G_q(\chi)=-1$ if $\chi$ is trivial. 

In general, explicit evaluations of Gauss sums are difficult. There are only a few cases that the Gauss sums have been completely evaluated. The most 
well-known case is  the {\it quadratic case}, i.e., the order of the multiplicative 
character involved is $2$. 
\begin{theorem}\label{thm:quad}{\em (\cite[Theorem~5.15]{LN97})}
Let $\eta$ be the quadratic character of $\F_{q}=\F_{p^s}$. Then,  
\[
G_{q}(\eta)
=(-1)^{s-1}\Big(\sqrt{(-1)^{\frac{p-1}{2}}p}\Big)^s. 
\]
\end{theorem}
The next simple case is the so-called {\it semi-primitive case}, where 
there exists an integer $\ell$ such that $p^\ell\equiv -1\,(\mod{N})$. Here, 
$N$ is the order of the multiplicative 
character involved. In particular, we give the following for 
later use. 
\begin{theorem}\label{thm:semi}{\em (\cite[Theorem~5.16]{LN97})}
Let $\chi$ be a nontrivial multiplicative character of $\F_{q^{2}}$ of order $N$ dividing $q+1$. Then, 
\begin{align*}
G_{q^{2}}(\chi)
=
\left\{
\begin{array}{ll}
q & \mbox{ if $N$ odd or $\tfrac{q+1}{N}$ even,}\\
-q, &\mbox{ if $N$ even and $\tfrac{q+1}{N}$ odd.}
 \end{array}
\right.
\end{align*}
\end{theorem}
We will also need the  {\it Davenport-Hasse product formula}, which is stated below. 
\begin{theorem}
\label{thm:Stickel2}{\em (\cite[Theorem~11.3.5]{BEW97})}
Let $\chi'$ be a multiplicative character of order $\ell>1$ of  $\F_{q}$. For  every nontrivial multiplicative character $\chi$ of $\F_{q}$, 
\[
G_q(\chi)=\frac{G_q(\chi^\ell)}{\chi^\ell(\ell)}
\prod_{i=1}^{\ell-1}
\frac{G_q({\chi'}^i)}{G_q(\chi{\chi'}^i)}. 
\]
\end{theorem}

Let $N$ be a positive integer dividing $q-1$. We set $C_i^{(N,q)}=\omega^i\langle \omega^N\rangle$, $0\leq i\leq N-1$, which are called the $N$th {\it cyclotomic classes} of $\F_q$. In this paper, we need to evaluate the (additive) character values of a union of some cyclotomic classes. 
In particular, the character sums defined by 
\[
\psi_{\F_q}(C_i^{(N,q)})=\sum_{x\in C_i^{(N,q)}}\psi_{\F_q}(x), \quad i=0,1,\ldots,N-1,
\]
are called the $N$th {\it Gauss periods} of $\F_q$. 
By the orthogonality of characters, the Gauss period can be expressed as a linear combination of Gauss sums: 
\begin{equation}\label{eq:ortho1}
\psi_{\F_{q}}(C_i^{(N,q)})=\frac{1}{N}\sum_{j=0}^{N-1}G_q(\chi^{j})\chi^{-j}(\omega^i), \, \quad i=0,1,\ldots,N-1,
\end{equation}
where $\chi$ is any fixed multiplicative character of order $N$ of $\F_q$. 
For example, if $N=2$, 
we have the following from Theorem~\ref{thm:quad}: 
\begin{equation}\label{eq:Gaussquad}
\psi_{\F_q}(C_i^{(2,q)})=\frac{-1+(-1)^iG_q(\eta)}{2}=\frac{-1+(-1)^{i+s-1+\frac{(p-1)s}{4}}p^\frac{s}{2}}{2}, \quad i=0,1, 
\end{equation}
where $\eta$ is the quadratic character of $\F_q$. On the other hand, the Gauss sum with respect to a multiplicative character $\chi$ of order $N$ can be expressed as a linear combination of Gauss periods: 
\begin{equation}\label{eq:ortho2}
G_{q}(\chi)=\sum_{i=0}^{N-1}\psi_{\F_{q}}(C_i^{(N,q)})\chi(\omega^i).  
\end{equation}
\subsection{Known results on projective sets of type Q}\label{sec:twonew}
Let $\PG(k-1,q)$ denote the $(k-1)$-dimensional projective space over $\F_q$. 
A set ${\mathcal S}$ of $n$ points of $\PG(k-1,q)$ is called a 
{\it projective $(n,k,h_1,h_2)$ set} if every hyperplane of $\PG(k-1,q)$ meets ${\mathcal S}$ in $h_1$ or $h_2$ points.  
In particular, a subset ${\mathcal S}$ of the point set of $\PG(3,q)$ is called {\it type Q} if 
\[
(n,k,h_1,h_2)=\Big(\frac{q^4-1}{4(q-1)},4,\frac{(q-1)^2}{4},\frac{(q+1)^2}{4}\Big).
\] 

In this paper, we will use the following model of $\PG(3,q)$: 
We view $\F_{q^2}\times \F_{q^2}$ as a $4$-dimensional vector space 
over $\F_q$. For a nonzero vector $(x,y)\in (\F_{q^2}\times \F_{q^2})\setminus \{(0,0)\}$, we use
$\langle (x,y)\rangle$ to denote the projective point in $\PG(3,q)$ corresponding to the one-dimensional subspace over $\F_q$ spanned by $(x,y)$. 
Let ${\mathcal P}$ be the set of points of $\PG(3,q)$.
Then, all (hyper)planes in $\PG(3,q)$ are given by 
\[
H_{a,b}=\{\langle (x,y)\rangle\,|\,\Tr_{q^2/q}(ax+by)=0\},\quad \langle(a,b)\rangle\in {\mathcal P}. 
\]
Let ${\mathcal S}$ be a set of points of $\PG(3,q)$, and  define 
\[
E=\{\lambda (x,y)\,|\,\lambda\in \F_{q}^\ast,\langle (x,y)\rangle\in {\mathcal S}\}. 
\]
Noting that each nontrivial additive character of $\F_{q^2}\times \F_{q^2}$ is given by  
\[
\psi_{a,b}((x,y))=\psi_{\F_{q^2}}(ax+by), \quad (x,y)\in \F_{q^2}\times \F_{q^2}, 
\]  
where $(0,0)\neq (a,b)\in  \F_{q^2}\times \F_{q^2}$, 
we have 
\begin{align*}
\psi_{a,b}(E)=&\,\sum_{\lambda \in \F_{q}}\sum_{\langle(x,y)\rangle \in {\mathcal S}}\psi_{\F_q}(\lambda\Tr_{q^2/q}(ax+by))-|{\mathcal S}|\\
=&\,q|H_{a,b}\cap {\mathcal S}|-|{\mathcal S}|. 
\end{align*}
Hence, we have the following proposition. 
\begin{proposition}\label{prop:twoint}
The set ${\mathcal S}$ is a projective set of type Q in $\PG(3,q)$ if and only if $|E|=\frac{q^4-1}{4}$ and 
$\psi_{a,b}(E)$ take exactly two values  $\frac{q^2-1}{4}$ and $\frac{-3q^2-1}{4}$ for all $(0,0)\neq (a,b)\in  \F_{q^2}\times \F_{q^2}$. 
\end{proposition}   
The set $E\subseteq \F_{q^2}\times \F_{q^2}$ is also called {\it type Q} if it satisfies the condition of Proposition~\ref{prop:twoint}.

\vspace{0.1in}

A {\it spread} in $\PG(3,q)$ is a collection ${\mathcal L}$ of $q^2+1$ pairwise 
skew lines; equivalently, ${\mathcal L}$ can be regarded as a 
collection ${\mathcal K}$ of 
$2$-dimensional subspaces of the underlying $4$-dimensional vector space 
$V(4,q)$ over $\F_q$, any two of which intersect at zero only. We also call such a set ${\mathcal K}$ of $2$-dimensional subspaces as a {\it spread} of $V(4,q)$. 

The following important theorem was given by Wilson and Xiang~\cite{WX97}.

\begin{theorem}\label{thm:HDiff}
Let ${\mathcal L}=\{L_i\,|\,0\le i\le q^2\}$ be a spread of 
$\PG(3,q)$, and  assume the existence of four pairwise disjoint projective sets ${\mathcal S}_i$, $i=1,2,3,4$, of type Q in $\PG(3,q)$ such that 
${\mathcal S}_0\cup {\mathcal S}_2=\bigcup_{i=0}^{(q^2-1)/2}L_i$ and 
${\mathcal S}_1\cup {\mathcal S}_3=\bigcup_{i=(q^2+1)/2}^{q^2}L_i$. 
Then there exists a Menon-Hadamard difference set in $H\times G$, where $H$ is either group of order $4$ and $G$ is an elementary abelian group of order $q^4$. 
\end{theorem}
\begin{remark}\label{rem:HDiff}
From Proposition~\ref{prop:twoint} and Theorem~\ref{thm:HDiff}, in order to construct a Menon-Hadamard difference set in a group of order $4q^4$, we need to 
find four disjoint sets $C_i\subseteq (\F_{q^2}\times \F_{q^2})\setminus \{(0,0)\}$, $i=0,1,2,3$, of type Q and a suitable spread 
${\mathcal K}=\{K_i\,|\,0\le i\le q^2\}$ consisting of $2$-dimensional subspaces  
of $V(4,q)$
such that $C_0\cup C_2 \cup \{(0,0)\}=\bigcup_{i=0}^{(q^2-1)/2}K_i$ and  $C_1\cup C_3\cup  \{(0,0)\}=\bigcup_{i=(q^2+1)/2}^{q^2}K_i$. 
\end{remark}
 
We now review the construction of projective sets of type Q given by Chen~\cite{Ch97}. Let $\omega$ be a primitive element of $\F_{q^2}$. Furthermore, let 
\begin{align*}
X=\{x\in \F_{q^2}\,|\,\Tr_{q^2/q}(x)\in C_0^{(2,q)}\}, \quad 
X'=\{x\omega \,|\,\Tr_{q^2/q}(x)\in C_0^{(2,q)}\}. 
\end{align*}
Define 
\begin{align*}
&X_1=X\setminus (X\cap X'),\, \, X_2=X'\setminus (X\cap X'),\\
&X_3=X\cap X',\, \, X_4=\F_{q^2}\setminus  (X_1\cup X_2\cup X_3),
\end{align*}
and 
\begin{align*}
C_0&\,=\{(x,xy)\,|\,x\in C_0^{(2,q^2)},y\in X_1\}\cup \{(x,xy)\,|\,x\in C_1^{(2,q^2)},y\in X_2\}\cup \{(0,x)\,|\,x\in C_\tau^{(2,q^2)}\},\\
C_1&\,=\{(x,xy)\,|\,x\in C_0^{(2,q^2)},y\in X_3\}\cup \{(x,xy)\,|\,x\in C_1^{(2,q^2)},y\in X_4\},\\
C_2&\,=\{(x,xy)\,|\,x\in C_1^{(2,q^2)},y\in X_1\}\cup \{(x,xy)\,|\,x\in C_0^{(2,q^2)},y\in X_2\}\cup \{(0,x)\,|\,x\in C_{\tau+1}^{(2,q^2)}\},\\
C_3&\,=\{(x,xy)\,|\,x\in C_1^{(2,q^2)},y\in X_3\}\cup \{(x,xy)\,|\,x\in C_0^{(2,q^2)},y\in X_4\},
\end{align*}
where $\tau=0$ or $1$ depending on whether $q\equiv 1$ or $3\,(\mod{4})$. It is clear that these type Q sets admit the automorphism $T$. 
\begin{theorem}\label{thm:chen}
The sets $C_i$, $i=0,1,2,3$, are type Q. Furthermore, these sets satisfy 
the assumption of Remark~\ref{rem:HDiff} with respect to the spread ${\mathcal K}$
consisting of the following $2$-dimensional subspaces: 
\[
K_y=\{(x,xy)\,|\,x\in \F_{q^2}\}, \, y \in \F_{q^2}, \mbox{ \, and \, }
K_{\infty}=\{(0,x)\,|\,x\in \F_{q^2}\}. 
\]
\end{theorem}

On the other hand, Wilson and Xiang~\cite{WX97} constructed Menon-Hadamard difference sets of order $4q^4$ for $q=5,13,17$ using the following four type Q sets: 
\begin{align*}
C_i=&\{(0,y)\,|\,y\in C_{\tau_i}^{(2,q^2)}\} \cup \{(xy,xy^{-1}\omega^j)\,|\,x\in \F_q^\ast,y\in C_0^{(2,q^2)},j\in A_i\}\\
&\quad \cup  
\{(xy,xy^{-1}\omega^j)\,|\,x\in \F_q^\ast,y\in C_1^{(2,q^2)},j\in B_i\}, \quad i=0,2,\\
C_i=&\{(y,0)\,|\,y\in C_{\tau_i}^{(2,q^2)}\}\cup \{(xy,xy^{-1}\omega^j)\,|\,x\in \F_q^\ast,y\in C_0^{(2,q^2)},j\in A_i\}\\
&\quad \cup  \{(xy,xy^{-1}\omega^j)\,|\,x\in \F_q^\ast,y\in C_1^{(2,q^2)},j\in B_i\}, \quad i=1,3,
\end{align*} 
for some subsets $A_i,B_i$, $i=0,1,2,3$, of $\{0,1,\ldots,2q+1\}$,  and 
the spread ${\mathcal K}$
consisting of the following $2$-dimensional subspaces: 
\[
K_y=\{(x,yx^q)\,|\,x\in \F_{q^2}\}, \, y \in \F_{q^2}, \mbox{ \, and \, }
K_{\infty}=\{(0,x)\,|\,x\in \F_{q^2}\}. 
\]
It is clear that these type Q sets admit the automorphism $T'$. 

\section{A generalization of Chen's construction}\label{sec:Chen}
We first fix notation used in this section. 
Let $q=p^s$ be an odd prime power with $p$ a prime, and $m$ be a fixed positive integer satisfying $2m\,|\,(q+1)$. Then, there exists a minimal $\ell$ such that $2m\,|\,(p^\ell+1)$. Write $s=\ell t$ for some $t\ge 1$. Let $\omega$ be a primitive element of $\F_{q^2}$.  
Let $T_i$, $i=0,1$, be two arbitrary subsets of $\F_{q}$, 
and
\begin{equation}\label{eq:defS01}
S_0=\{x\,|\,\Tr_{q^2/q}(x)\in T_0\},\quad
S_1=\{x\,|\,\Tr_{q^2/q}(x\omega^m)\in T_1\}. 
\end{equation}
Furthermore, let $K$ be any $m$-subset of $\{0,1,\ldots,2m-1\}$ such that 
$K\cap \{x+m\,(\mod{2m})\,|\,x \in K\}=\emptyset$. 
Define 
\begin{equation}\label{eq:defA01}
A_0=S_0\setminus S_1, \quad A_1=S_1\setminus S_0,\quad 
D_0=\bigcup_{i\in K}C_i^{(2m,q^2)},\quad
D_1=\bigcup_{i\in K}C_{i+m}^{(2m,q^2)},
\end{equation}
and 
\[
\epsilon:=\left\{
\begin{array}{ll}
1, &\mbox{ if $(p^\ell+1)/2m$ is even and $t$ is odd,}\\
0, &\mbox{  otherwise. }
 \end{array}
\right.
\]
\begin{remark}\label{rem:secChen}
\begin{itemize}
\item[(i)] The indicator function of $S_i$, $i=0,1$, is given by 
\[
f_{S_i}(y)=\frac{1}{q}\sum_{c\in \F_q}\sum_{u\in T_i}\psi_{\F_{q^2}}(cy\omega^{mi})
\psi_{\F_q}(-cu), \quad i=0,1. 
\]
\item[(ii)] The size of each $S_i$ is $q|T_i|$ since $\Tr_{q^2/q}$ is a linear mapping over $\F_{q}$.
\item[(iii)] The size of $S_0\cap S_1$ is $|T_0||T_1|$; it is clear that 
\begin{align}
|S_0\cap S_1|=&\,\sum_{y\in \F_{q^2}}f_{S_0}(y)f_{S_1}(y)\nonumber\\
=&\,\frac{1}{q^2}\sum_{c,d\in \F_q}\sum_{u\in T_0}\sum_{v\in T_1}\sum_{y\in \F_{q^2}}\psi_{\F_{q^2}}(y(c+d\omega^{m}))
\psi_{\F_q}(-cu-dv). \label{eq:sizes0s1}
\end{align}
Since $\omega^m\not\in \F_{q}$, $c+d\omega^m=0$ if and only if $c=d=0$. Hence, the right-hand side of \eqref{eq:sizes0s1} is equal to $|T_0||T_1|$. 
\item[(iv)] 
Since $2m\,|\,(q+1)$,  the character values of $D_i\subseteq \F_{q^2}$, $i=0,1$, can be evaluated by using \eqref{eq:ortho1} and the Gauss sums in semi-primitive case (see, e.g.,  \cite[Theorem~2]{bwx}): 
for $b\in \F_{q^2}^\ast$, 
\[
\sum_{x\in D_\epsilon}\psi_{\F_{q^2}}(bx)=\left\{
\begin{array}{ll}
\frac{-1-q}{2}, & \mbox{ if $b^{-1}\in D_0$,}\\
\frac{-1+q}{2}, &\mbox{ if  $b^{-1}\in D_1$. }
 \end{array}
\right. 
\]
\end{itemize}
\end{remark}
The following is our main result in this section. 
\begin{theorem}\label{thm:mainf1}
\begin{itemize}
\item[(1)] Assume that $|T_0|=|T_1|=(q-1)/2$, and define 
\[
E_0=\{(x,xy)\,|\,x\in D_0,y \in A_0\} \cup  
\{(x,xy)\,|\,x\in D_1,y \in A_1\} \cup \{(0,x)\,|\,x \in D_\epsilon\}. 
\] 
Then $E_0$ is a set of type Q in $\F_{q^2}\times \F_{q^2}$. 
\item[(2)] Assume that $|T_0|=(q-1)/2$ and $|T_1|=(q+1)/2$, 
and define 
\[
E_1=\{(x,xy)\,|\,x\in D_0,y \in A_0\} \cup  
\{(x,xy)\,|\,x\in D_1,y \in A_1\}. 
\] 
Then $E_1$ is a set of type Q in $\F_{q^2}\times \F_{q^2}$. 
\end{itemize}
\end{theorem}
This theorem obviously generalizes the construction of type Q sets given 
by Chen~\cite{Ch97}. Indeed, we used $D_i$, $i=0,1$, instead of $C_i^{(2,q^2)}$, $i=0,1$, in the definition of $X$ and $X'$ (see Subsection~\ref{sec:twonew}). This new construction is much more flexible than that in \cite{Ch97}. 

To prove this theorem, 
we will evaluate the character values $\psi_{a,b}(E_i)$, $(a,b)\in (\F_{q^2}\times \F_{q^2})\setminus \{(0,0)\}$, by a series of the following lemmas. We first treat the case where $b=0$. 
\begin{lemma}\label{lem:b0}
For $b=0$ and $a\not=0$, it holds that 
\[
\psi_{a,b}(E_0)=\frac{q^2-1}{4}. 
\]
\end{lemma}
\proof
Since $|T_0|=|T_1|=(q-1)/2$, by Remark~\ref{rem:secChen}~(ii),(iii), 
we have $|A_0|=|A_1|=(q^2-1)/4$. Then, we have 
\begin{align*}
\psi_{a,0}(E_0)=&\,\sum_{x\in D_0}\sum_{y\in A_0}\psi_{\F_{q^2}}(ax)+
\sum_{x\in D_1}\sum_{y\in A_1}\psi_{\F_{q^2}}(ax)+\frac{q^2-1}{2}\\
=&\frac{q^2-1}{4}\sum_{x\in \F_{q^2}^\ast}\psi_{\F_{q^2}}(ax)+\frac{q^2-1}{2}=
\frac{q^2-1}{4}. 
\end{align*}
This completes the proof. 
\qed
\begin{lemma}\label{lem:b02}
For $b=0$ and $a\not=0$, we have 
\[
\psi_{a,b}(E_1)=\left\{
\begin{array}{ll}
\frac{q^2-1}{4}, & \mbox{ if $a^{-1}\in D_\epsilon$,}\\
\frac{-3q^2-1}{4}, &\mbox{ otherwise.}
 \end{array}
\right. 
\]
\end{lemma}
\proof
Since $|T_0|=(q-1)/2$ and $|T_1|=(q+1)/2$, by Remark~\ref{rem:secChen}~(ii),(iii),  we have $|A_0|=(q-1)^2/4$ and $|A_1|=(q+1)^2/4$. Then, we have 
\begin{align}
\psi_{a,0}(E_1)=&\,\sum_{x\in D_0}\sum_{y\in A_0}\psi_{\F_{q^2}}(ax)+
\sum_{x\in D_1}\sum_{y\in A_1}\psi_{\F_{q^2}}(ax)\nonumber\\
=&\frac{(q-1)^2}{4}\sum_{x\in \F_{q^2}^\ast}\psi_{\F_{q^2}}(ax)+q\sum_{x\in D_1}\psi_{\F_{q^2}}(ax). \label{eq:b0ane0}
\end{align}
Finally, by Remark~\ref{rem:secChen} (iv), \eqref{eq:b0ane0} is reformulated as 
\[
\psi_{a,0}(E_1)=-\frac{(q-1)^2}{4}+q
\left\{
\begin{array}{ll}
\frac{-1+q}{2}, & \mbox{ if $a^{-1}\in D_\epsilon$,}\\
\frac{-1-q}{2}, &\mbox{ otherwise.}
 \end{array}
\right.
\]
This completes the proof. 
\qed
\vspace{0.3cm}

We next treat the case where $b\not=0$. 
Let $f_{S_i}$, $i=0,1$, be defined as in Remark~\ref{rem:secChen}~(i). Define 
\begin{align*}
U_1=&\sum_{x\in D_0}\sum_{y\in \F_{q^2}}\psi_{\F_{q^2}}(x(a+by))
f_{S_0}(y),\\
U_2=&\sum_{x\in D_1}\sum_{y\in \F_{q^2}}\psi_{\F_{q^2}}(x(a+by))
f_{S_1}(y),\\
U_3=&\sum_{x\in \F_{q^2}^\ast}\sum_{y\in \F_{q^2}}\psi_{\F_{q^2}}(x(a+by))
f_{S_0}(y)f_{S_1}(y). 
\end{align*}
Then, the character values of $E_i$, $i=0,1$, are given by  
\begin{equation}\label{eq:chraE0}
\psi_{a,b}(E_0)=
U_1+U_2-U_3+\sum_{x\in D_\epsilon}\psi_{\F_{q^2}}(bx)
\end{equation}
and 
\begin{equation}\label{eq:chraE1}
\psi_{a,b}(E_1)=
U_1+U_2-U_3. 
\end{equation}

\begin{lemma}\label{lem:u1}
If $b\not=0$, it holds that 
\[
U_1=
\left\{
\begin{array}{ll}
-q|T_0|+q^2, &\mbox{ if $-ab^{-1}\in S_0$ and $b^{-1}\in D_0$,}\\
-q|T_0|, &\mbox{  if $-ab^{-1}\not\in S_0$ and $b^{-1}\in D_0$,}\\
0, &\mbox{  if $b^{-1}\in D_1$.}\\
 \end{array}
\right.
\]
\end{lemma}
\proof 
If $b\not=0$, we have 
\begin{equation}\label{eq:U1}
U_1=\frac{1}{q}\sum_{x\in D_0}\sum_{y\in \F_{q^2}}
\sum_{c\in \F_q}\sum_{u\in T_0}\psi_{\F_{q^2}}(xa)
\psi_{\F_{q^2}}((xb+c)y)
\psi_{\F_q}(-cu). 
\end{equation}
If $b^{-1}\in D_1$, there are no $x\in D_0$ such that $xb+c=0$; we have 
$U_1=0$. If $b^{-1} \in D_0$, continuing from \eqref{eq:U1}, we have 
\begin{align*}
U_1=&\,q
\sum_{c\in \F_q^\ast}\sum_{u\in T_0}\psi_{\F_{q^2}}(-acb^{-1})
\psi_{\F_q}(-cu)\\
=&\,-q|T_0|+q
\sum_{c\in \F_q}\sum_{u\in T_0}\psi_{\F_q}(\Tr_{q^2/q}(-ab^{-1})c-cu)\\
=&\,-q|T_0|+q^2
\left\{
\begin{array}{ll}
1, & \mbox{ if $\Tr_{q^2/q}(-ab^{-1})\in T_0$,}\\
0, &\mbox{ otherwise.}
 \end{array}
\right.
\end{align*}
This completes the proof. 
\qed
\vspace{0.3cm}

\begin{lemma}\label{lem:u2}
If $b\not=0$, we have
\[
U_2=
\left\{
\begin{array}{ll}
-q|T_1|+q^2, &\mbox{ if $-ab^{-1}\in S_1$ and $b^{-1}\in D_0$,}\\
-q|T_1|, &\mbox{ if $-ab^{-1}\not\in S_1$ and $b^{-1}\in D_0$,}\\
0, &\mbox{ if $b^{-1}\in D_1$.}
 \end{array}
\right.
\]
\end{lemma}
\proof 
If $b\not=0$, we have 
\begin{equation}\label{eq:U2}
U_2=\frac{1}{q}\sum_{x\in D_1}\sum_{y\in \F_{q^2}}
\sum_{c\in \F_q}\sum_{u\in T_1}\psi_{\F_{q^2}}(xa)
\psi_{\F_{q^2}}((xb+c\omega^{m})y)
\psi_{\F_q}(-cu). 
\end{equation}
If $b^{-1}\in D_1$, there are no $x\in D_1$ such that $xb+c\omega^{m}=0$; hence 
$U_2=0$. If $b^{-1}\in D_0$, continuing from \eqref{eq:U2}, we have 
\begin{align*}
U_2=&\,q
\sum_{c\in \F_q^\ast}\sum_{u\in T_1}\psi_{\F_{q^2}}(-acb^{-1}\omega^{m})
\psi_{\F_q}(-cu)\\
=&\,-q|T_1|+q
\sum_{c\in \F_q}\sum_{u\in T_1}\psi_{\F_q}(\Tr_{q^2/q}(-ab^{-1}\omega^{m})c-cu)\\
=&\,-q|T_1|+q^2
\left\{
\begin{array}{ll}
1, & \mbox{ if $\Tr_{q^2/q}(-ab^{-1}\omega^{m})\in T_1$,}\\
0, &\mbox{ otherwise.}
 \end{array}
\right.
\end{align*}
This completes the proof. 
\qed
\begin{lemma}\label{lem:u3}
If $b\not=0$, we have 
\[
U_3=
\left\{
\begin{array}{ll}
-|T_0||T_1|+q^2, & \mbox{ if $-ab^{-1}\in S_0\cap S_1$,}\\
-|T_0||T_1|, &\mbox{ otherwise.}
 \end{array}
\right.
\]
\end{lemma}
\proof 
Note that $D_0\cup D_1=\F_{q^2}^\ast$ and $|S_0\cap S_1|=|T_0||T_1|$. 
Since $b\not=0$, we have 
\begin{align*}
U_3=&\,\sum_{x\in \F_{q^2}}\sum_{y\in \F_{q^2}}\psi_{\F_{q^2}}(x(a+by))
f_{S_0}(y)f_{S_1}(y)-|S_0\cap S_1|\\
=&\,q^2
f_{S_0}(-ab^{-1})f_{S_1}(-ab^{-1})-|T_0||T_1|. 
\end{align*}
This completes the proof.  
\qed
\vspace{0.3cm}

{\bf Proof of Theorem~\ref{thm:mainf1}:}\, 
In the case where $b=0$, the statement follows from Lemmas~\ref{lem:b0} and \ref{lem:b02}. We now treat the case where $b\not=0$. 
By the evaluations for $U_1,U_2,U_3$ in Lemmas~\ref{lem:u1}--\ref{lem:u3}, 
we have 
\begin{align*}
U_1+U_2-U_3=
\left\{
\begin{array}{ll}
-q(|T_0|+|T_1|-q)+|T_0||T_1|, & \mbox{if $b^{-1}\in D_0$, $-ab^{-1}\in S_0$,  $-ab^{-1}\in S_1$;}\\
 & \mbox{\, or $b^{-1}\in D_0$, $-ab^{-1}\not\in S_0$,  $-ab^{-1}\in S_1$;}\\
 & \mbox{\, or $b^{-1}\in D_0$, $-ab^{-1}\in S_0$, $-ab^{-1}\not\in S_1$,}\\
-q(|T_0|+|T_1|)+|T_0||T_1|, &\mbox{if  $b^{-1}\in D_0$, $-ab^{-1}\not\in S_0$,  $-ab^{-1}\not\in S_1$, }\\
-q^2+|T_0||T_1|, & \mbox{if $b^{-1}\in D_1$, $-ab^{-1}\in S_0$,  $-ab^{-1}\in S_1$,}\\
|T_0||T_1|, &\mbox{if  $b^{-1}\in D_1$, $-ab^{-1}\not\in S_0$,  $-ab^{-1}\in S_1$;}\\
 &\mbox{\, or  $b^{-1}\in D_1$, $-ab^{-1}\in S_0$,   $-ab^{-1}\not\in S_1$;}\\
 &\mbox{\, or  $b^{-1}\in D_1$, $-ab^{-1}\not\in S_0$,  $-ab^{-1}\not\in S_1$.}
 \end{array}
\right. 
\end{align*} 
(1) Since $|T_0|=|T_1|=(q-1)/2$, by Remark~\ref{rem:secChen} (iv), we have
\begin{align*}
\psi_{a,b}(E_0)=&\,U_1+U_2-U_3+\sum_{x\in D_\epsilon}\psi_{\F_{q^2}}(bx)\\
=&\,
\left\{
\begin{array}{ll}
\frac{-3q^2-1}{4}, & \mbox{ if  $b^{-1}\in D_0$, $-ab^{-1}\not\in S_0$, and  $-ab^{-1}\not\in S_1$;}\\
 & \mbox{ or  if $b^{-1}\in D_1$, $-ab^{-1}\in S_0$, and  $-ab^{-1}\in S_1$,}\\
\frac{q^2-1}{4}, & \mbox{ otherwise.}
 \end{array}
\right. 
\end{align*}   
(2) Since $|T_0|=(q-1)/2$ and $|T_1|=(q+1)/2$, we have
\begin{align*}
\psi_{a,b}(E_1)=&\,U_1+U_2-U_3\\
=&\,
\left\{
\begin{array}{ll}
\frac{-3q^2-1}{4}, & \mbox{ if  $b^{-1}\in D_0$, $-ab^{-1}\not\in S_0$, and  $-ab^{-1}\not\in S_1$;}\\
 & \mbox{ or  if $b^{-1}\in D_1$, $-ab^{-1}\in S_0$, and  $-ab^{-1}\in S_1$,}\\
\frac{q^2-1}{4}, & \mbox{ otherwise.}
 \end{array}
\right. 
\end{align*}   
This completes the proof of the theorem. \qed 
\begin{corollary}
Let $T_i$, $i=0,1$, be arbitrary $(q-1)/2$-subsets of $\F_q$ and $S_0,S_1,A_0,A_1$ be the sets defined as in \eqref{eq:defS01} and \eqref{eq:defA01}. 
Furthermore, define 
\[
S_1'=\{x\in \F_{q^2}\,|\,\Tr_{q^2/q}(x\omega^m)\in \F_q\setminus T_1\}, \quad  
A_0'=S_0\setminus S_1', \quad A_1'=S_1'\setminus S_0. 
\]
Then, 
the sets 
\begin{align*}
C_0=&\{(x,xy)\,|\,x\in D_0,y \in A_0\} \cup  
\{(x,xy)\,|\,x\in D_1,y \in A_1\} \cup \{(0,x)\,|\,x \in D_\epsilon\},\\ 
C_1=&\{(x,xy)\,|\,x\in D_0,y \in A_0'\} \cup  
\{(x,xy)\,|\,x\in D_1,y \in A_1'\},\\ 
C_2=&\{(x,xy)\,|\,x\in D_1,y \in A_0\} \cup  
\{(x,xy)\,|\,x\in D_0,y \in A_1\} \cup \{(0,x)\,|\,x \in D_{\epsilon+1}\},\\ 
C_3=&\{(x,xy)\,|\,x\in D_1,y \in A_0'\} \cup  
\{(x,xy)\,|\,x\in D_0,y \in A_1'\} 
\end{align*} 
are of type Q, where 
the subscript of $D_{\epsilon+1}$ is reduced modulo $2$.  Furthermore, these sets satisfy the assumptions of Remark~\ref{rem:HDiff} with respect to the spread consisting of the following $2$-dimensional subspaces: 
\[
K_y=\{(x,xy)\,|\,x\in \F_{q^2}\}, y\in \F_{q^2}, \mbox{\, and \, }K_\infty=\{(0,x)\,|\,x\in \F_{q^2}\}.  
\]
\end{corollary}
\proof 
By 
Theorem~\ref{thm:mainf1}, $C_0$ and $C_1$ are type Q sets. Furthermore, since $C_2=\omega^m C_0$ and $C_3=\omega^m C_1$, the sets $C_2$ and $C_3$ are also  of type Q. Finally,  $C_i$, $i=0,1,2,3$, satisfy 
the assumption of Remark~\ref{rem:HDiff} as $C_0\cup C_2\cup\{(0,0)\}=\big(\bigcup_{y\in A_0\cup A_1}K_y\big)\cup K_\infty$ and $C_1\cup C_3\cup\{(0,0)\}=\bigcup_{y\in A_0'\cup A_1'}K_y$. 
\qed

\section{A generalization of Wilson-Xiang's examples}\label{sec:firstty}
\subsection{The Setting}\label{subsec:set}
We fix notation used in this section. 
Let $q$ be a prime power  
and 
$\omega$ be a  primitive element of $\F_{q^2}$.  
Let $c$ be any fixed odd integer in $\{0,1,\ldots,2q+1\}=\Z_{2q+2}$.  

Define the following subsets of $\{0,1,\ldots,2q+1\}$:  
\begin{align*}
I_1=&\,\{i\,(\mod{2(q+1)})\,|\,\Tr_{q^2/q}(\omega^i)=0\}=\{\tfrac{q+1}{2},\tfrac{3(q+1)}{2}\},\\
I_2=&\,\{i\,(\mod{2(q+1)})\,|\,\Tr_{q^2/q}(\omega^i)\in C_0^{(2,q)}\},\\
I_3=&\,\{i\,(\mod{2(q+1)})\,|\,\Tr_{q^2/q}(\omega^i)\in  C_1^{(2,q)}\},\\
J_{i}=&\,I_i-c\,(\mod{2(q+1)}),\, \, \, \, i=1,2,3.  
\end{align*}
Then $|I_1|=2$, $|I_2|=|I_3|=q$, and $I_1\cup I_2\cup I_3=\Z_{2q+2}$.
Furthermore, define 
\begin{align*}
X_{1,c}=&\,(I_1\cap J_{2})\cup (I_2\cap J_{1}),\\
X_{2,c}=&\,(I_1\cap J_{3})\cup (I_3\cap J_{1}),\\
X_{3,c}=&\,I_2\cap J_{2},\, \, X_{4,c}=I_3\cap J_{3},\\
X_{5,c}=&\,(I_2\cap J_{3})\cup (I_3\cap J_{2}).
\end{align*}
It is clear that the $X_{i,c}$'s partition $\Z_{2q+2}$.  
In the appendix, we will show that the $X_{i,c}$'s have the following properties: 
\begin{itemize}
\item[(P1)] $X_{1,c}\equiv X_{2,c}+(q+1)\,(\mod{2(q+1)})$, $X_{3,c}\equiv X_{4,c}+(q+1)\,(\mod{2(q+1)})$; 
\item[(P2)] $|X_{1,c}|=|X_{2,c}|=2$, $|X_{3,c}|=|X_{4,c}|=\frac{q-1}{2}$, $|X_{5,c}|=q-1$; 
\item[(P3)] $X_{3,c+q+1}\cup X_{4,c+q+1}=X_{5,c}$;
\item[(P4)] $X_{1,c}+c\equiv -X_{1,c}+(q+1)\,(\mod{2(q+1)})$ or 
$X_{2,c}+c\equiv -X_{2,c}+(q+1)\,(\mod{2(q+1)})$ according as $q\equiv 3$ or $1\,(\mod{4})$. 
\item[(P5)] $|X_{1,c}\cap X_{1,c+q+1}|=1$; 
\item[(P6)] By the properties~(P2) and (P5), we can assume that $X_{1,c}=\{\alpha,\beta\}$ and $X_{1,c+q+1}=\{\alpha,\gamma\}$. Then, $\beta\equiv \gamma+(q+1)\,(\mod{2(q+1)})$. Furthermore, $\alpha\equiv 0\,(\mod{2})$ and $\beta\equiv 1\,(\mod{2})$ or  $\alpha\equiv 1\,(\mod{2})$ and $\beta\equiv 0\,(\mod{2})$ according as $q\equiv 3$ or $1\,(\mod{4})$. 
\item[(P7)] Define $R_i=\bigcup_{j\in X_{i,c}}C_{j}^{(2(q+1),q^2)}$, $i=1,2,3,4,5$. Then, 
$R_i$ takes the character values listed in Table~\ref{tab_1}: 
{\footnotesize
\begin{table}[!h]
\begin{center}
\caption{\label{tab_1}The values of $\psi_{\F_{q^2}}(\omega^a R_i)$'s}
\begin{tabular}{|c|c|c|c|c|c|}\hline
   & $R_1$& $R_2$ &$R_3$&$R_4$&$R_5$\\ \hline
$a\in Y_{1,c}$& $\frac{-2+q+G_q(\eta)}{2}$ & $\frac{-2+q- G_q(\eta)}{2}$ &  $\frac{(q-1)(-1+ G_q(\eta))}{4}$ & $\frac{(q-1)(-1- G_q(\eta))}{4}$ & $\frac{-q+1}{2}$\\\hline
 $a\in Y_{2,c}$ & $\frac{-2+q- G_q(\eta)}{2}$& $\frac{-2+q+ G_q(\eta)}{2}$ &  $\frac{(q-1)(-1- G_q(\eta))}{4}$ &  $\frac{(q-1)(-1+ G_q(\eta))}{4}$& $\frac{-q+1}{2}$ \\\hline
 $a\in Y_{3,c}$& $-1+ G_q(\eta)$ &  $-1- G_q(\eta)$ & $\frac{(1- G_q(\eta))^2}{4}$ & $\frac{(1+  G_q(\eta))^2}{4}$& $\frac{1-(-1)^{\frac{q-1}{2}}q}{2}$\\\hline
 $a\in Y_{4,c}$ &  $-1- G_q(\eta)$ &  $-1+ G_q(\eta)$ &  $\frac{(1+  G_q(\eta))^2}{4}$ &   $\frac{(1- G_q(\eta))^2}{4}$& 
$\frac{1-(-1)^{\frac{q-1}{2}}q}{2}$\\\hline
 $a\in Y_{5,c}$& $-1$ & $-1$ & $\frac{1-(-1)^{\frac{q-1}{2}}q}{4}$ & $\frac{1-(-1)^{\frac{q-1}{2}}q}{4}$& 
$\frac{1+(-1)^{\frac{q-1}{2}}q}{2}$\\\hline
\end{tabular}
\end{center}
\end{table}}
In the language of association schemes, the Cayley graphs on $(\F_{q^2},+)$ with connection sets $R_i$'s, together with the diagonal relation arising from the connection set $R_0=\{0\}$, form a $5$-class translation association scheme. 
Here, $Y_{i,c}$'s are subsets of $\{0,1,\ldots,2q+1\}$ defined as the 
index sets of the dual association scheme. 
\item[(P8)] $-Y_{i,c}+c\equiv Y_{i,c}\,(\mod{2(q+1)})$, $i=1,2$;
\item[(P9)] $-(Y_{3,c}\cup Y_{4,c})\equiv Y_{5,c}-c\,(\mod{2(q+1)})$;
\item[(P10)] Define $R_i'=\bigcup_{j\in Y_{i,c}}C_{j}^{(2(q+1),q^2)}$, $i=1,2,3,4,5$. Then, 
$R_i'$ takes the character values listed in Table~\ref{tab_2}: 
{\footnotesize
\begin{table}[!h]
\begin{center}
\caption{\label{tab_2}The values of $\psi_{\F_{q^2}}(\omega^a R_i')$'s}
\begin{tabular}{|c|c|c|c|c|c|}\hline
   & $R_1'$& $R_2'$ &$R_3'$&$R_4'$&$R_5'$\\ \hline
$a\in X_{1,c}$& $\frac{-2+q+G_q(\eta)}{2}$ & $\frac{-2+q- G_q(\eta)}{2}$ &  $\frac{(q-1)(-1+ G_q(\eta))}{4}$ & $\frac{(q-1)(-1- G_q(\eta))}{4}$ & $\frac{-q+1}{2}$\\\hline
 $a\in X_{2,c}$ & $\frac{-2+q- G_q(\eta)}{2}$& $\frac{-2+q+ G_q(\eta)}{2}$ &  $\frac{(q-1)(-1- G_q(\eta))}{4}$ &  $\frac{(q-1)(-1+ G_q(\eta))}{4}$& $\frac{-q+1}{2}$ \\\hline
 $a\in X_{3,c}$& $-1+ G_q(\eta)$ &  $-1- G_q(\eta)$ & $\frac{(1- G_q(\eta))^2}{4}$ & $\frac{(1+  G_q(\eta))^2}{4}$& $\frac{1-(-1)^{\frac{q-1}{2}}q}{2}$\\\hline
 $a\in X_{4,c}$ &  $-1- G_q(\eta)$ &  $-1+ G_q(\eta)$ &  $\frac{(1+  G_q(\eta))^2}{4}$ &   $\frac{(1- G_q(\eta))^2}{4}$& 
$\frac{1-(-1)^{\frac{q-1}{2}}q}{2}$\\\hline
 $a\in X_{5,c}$& $-1$ & $-1$ & $\frac{1-(-1)^{\frac{q-1}{2}}q}{4}$ & $\frac{1-(-1)^{\frac{q-1}{2}}q}{4}$& 
$\frac{1+(-1)^{\frac{q-1}{2}}q}{2}$\\\hline
\end{tabular}
\end{center}
\end{table}}
\end{itemize}
\subsection{The Construction}\label{subsec:const}
Let $X_{i,c}$, $Y_{i,c}$, $R_i$, $R_i'$, $i=1,2,\ldots,5$, be sets defined as in 
Subsection~\ref{subsec:set}. 
Let $A$ and $B$ be subsets of $\{0,1,\ldots,2q+1\}$ satisfying 
$A\cap B=X_{3,c}$ and as multisets, $A\cup B=X_{1,c}\cup X_{3,c}\cup X_{3,c}$. It follows that $(A\setminus B)\cup (B\setminus A)=X_{1,c}$.

Let $\tau=0$ or $1$ according as $q\equiv 3$ or $1\,(\mod{4})$. Define 
\begin{align}
D_0=&\{(0,y)\,|\,y\in C_\tau^{(2,q^2)}\},\nonumber\\
D_1=&\{(y,0)\,|\,y\in C_0^{(2,q^2)}\},\nonumber\\
D_2=&\{(xy,xy^{-1}\omega^i)\,|\,x\in \F_q^\ast,y\in C_0^{(2,q^2)},i\in A \},
\nonumber\\
D_3=&
\{(xy,xy^{-1}\omega^i)\,|\,x\in \F_q^\ast,y\in C_1^{(2,q^2)},i\in B \}.  \label{def:Dis}
\end{align}

We denote the set of even (resp. odd) elements in any subset $S$ of $\{0,1,\ldots,2q+1\}$ by $S_e$ (resp. $S_o$). 
The following is our main result in this section. 
\begin{theorem} \label{thm:mainWX1}
\begin{itemize}
\item[(1)] If $|A|=|B|=\frac{q+1}{2}$ and $|A_e|+|B_o|=|A_o|+|B_e|-2(-1)^{\frac{q-1}{2}}$, then $E_0=D_0\cup D_2\cup D_3$ 
is a type Q set in $\F_{q^2}\times \F_{q^2}$. 
\item[(2)] If $|A|=\frac{q+3}{2}$, $|B|=\frac{q-1}{2}$ and $|A_e|+|B_o|=|A_o|+|B_e|$, then  $E_1=D_1\cup D_2\cup D_3$ 
is a type Q set in $\F_{q^2}\times \F_{q^2}$. 
\end{itemize}
\end{theorem}
This theorem  generalizes the examples of type Q sets found by Wilson-Xiang~\cite{WX97}. Indeed, these sets admit the automorphism $T'$. See Subsection~\ref{sec:twonew}. 

To prove the theorem above, we will evaluate the character values of $E_i$, $i=0,1$. 
Define 
\begin{align*}
V_0=&\sum_{y\in C_\tau^{(2,q^2)}}\psi_{\F_{q^2}}(by), \quad 
V_1=\sum_{y\in C_0^{(2,q^2)}}\psi_{\F_{q^2}}(ay), \\
V_2=&\frac{1}{2}\sum_{i\in A}\sum_{y\in C_0^{(2,q^2)}}\sum_{x\in \F_q^\ast}\psi_{\F_{q^2}}(axy)\psi_{\F_{q^2}}(bxy^{-1}\omega^i), \\
V_3=&\frac{1}{2}\sum_{i\in B}\sum_{y\in C_1^{(2,q^2)}}\sum_{x\in \F_q^\ast}\psi_{\F_{q^2}}(axy)\psi_{\F_{q^2}}(bxy^{-1}\omega^i).
\end{align*}
Noting that each element in $D_2$ (resp. $D_3$) appears exactly twice when $x$ runs through $\F_q^*$ and $y$ runs through $C_0^{(2,q^2)}$ (resp. $C_1^{(2,q^2)}$),  we have $\psi_{a,b}(E_0)=V_0+V_2+V_3$ and $\psi_{a,b}(E_1)=V_1+V_2+V_3$. 
We will evaluate these character sums by considering two cases: (i) exactly one of $a, b$ is zero; and (ii) $a\not=0$ and $b\not=0$.  We first treat Case (i). 
\begin{lemma}\label{lemma:ab0}
If exactly one of $a, b$ is zero, then 
\[
\psi_{a,b}(E_0)=\left\{
\begin{array}{ll}
\frac{-3q^2-1}{4}, & \mbox{ if  $a=0$ and $b\in C_{1}^{(2,q^2)}$, }\\
\frac{q^2-1}{4}, & \mbox{ otherwise.}
 \end{array}
\right. 
\]
\end{lemma}
\proof 
If  $a\not=0$ and  $b=0$,  
it is clear that $V_0=\frac{q^2-1}{2}$. 
Furthermore, since $|A|=|B|=\frac{q+1}{2}$, we have 
\[
V_2+V_3=\frac{q+1}{4}\sum_{y\in \F_{q^2}^\ast}\sum_{x\in \F_q^\ast}\psi_{\F_{q^2}}(axy)=-\frac{q^2-1}{4}.
\]
Hence, $\psi_{a,b}(E_0)=\frac{q^2-1}{4}$. 
If $a=0$ and $b\not=0$, we have
\begin{align}\label{eq:v1v2v3}
&V_0+V_2+V_3\nonumber\\
=&\, 
\psi_{\F_{q^2}}(b C_\tau^{(2,q^2)})+\frac{q-1}{2}((|A_e|+|B_o|)
\psi_{\F_{q^2}}(b C_0^{(2,q^2)})+(|A_o|+ |B_e|)\psi_{\F_{q^2}}(b C_1^{(2,q^2)})). 
\end{align}
Since $|A_e|+|A_o|+|B_e|+|B_o|=|A|+|B|=q+1$ and $|A_e|+|B_o|=|A_o|+|B_e|-2(-1)^{\frac{q-1}{2}}$, we have $|A_e|+|B_o|=\frac{q+1}{2}-(-1)^{\frac{q-1}{2}}$ and $|A_o|+|B_e|=\frac{q+1}{2}+(-1)^{\frac{q-1}{2}}$. 
Hence, \eqref{eq:v1v2v3} is reformulated as 
\[
V_0+V_2+V_3=q\psi_{\F_{q^2}}(b C_\tau^{(2,q^2)})-\Big(\frac{q-1}{2}\Big)^2.
\]
Finally, by \eqref{eq:Gaussquad}, the statement follows. 
\qed
\begin{lemma}\label{lemma:ab02}
If exactly one of $a, b$ is zero, then 
\[
\psi_{a,b}(E_1)=\left\{
\begin{array}{ll}
\frac{-3q^2-1}{4}, & \mbox{ if  $b=0$ and $a\in C_{\tau+1}^{(2,q^2)}$, }\\
\frac{q^2-1}{4}, & \mbox{ otherwise.}
 \end{array}
\right. 
\]
\end{lemma}
\proof 
If $a=0$ and $b\not=0$, it is clear that 
$V_1=\frac{q^2-1}{2}$. Since $|A_e|+|B_o|=|A_o|+|B_e|=\frac{q+1}{2}$,  we have
\begin{align*}
V_2+V_3
=&\,\frac{q-1}{2}((|A_e|+|B_o|)
\psi_{\F_{q^2}}(b C_0^{(2,q^2)})+(|A_o|+ |B_e|)\psi_{\F_{q^2}}(b C_1^{(2,q^2)})) \\
=&\,-\frac{q^2-1}{4}. 
\end{align*}
Hence, $\psi_{a,b}(E_1)=\frac{q^2-1}{4}$. 
If  $a\not=0$ and  $b=0$,  
\begin{equation}\label{eq:v1v2v3R}
V_1+V_2+V_3=\psi_{\F_{q^2}}(aC_0^{(2,q^2)})
+\frac{q-1}{2}(|A|\psi_{\F_{q^2}}(aC_0^{(2,q^2)})+|B|\psi_{\F_{q^2}}(aC_1^{(2,q^2)})). 
\end{equation}
Since $|A|=\frac{q+3}{2}$ and $|B|=\frac{q-1}{2}$, \eqref{eq:v1v2v3R} is reformulated as 
\[
V_1+V_2+V_3=q\psi_{\F_{q^2}}(a C_0^{(2,q^2)})-\Big(\frac{q-1}{2}\Big)^2. 
\]
Finally, by \eqref{eq:Gaussquad}, the statement follows. 
\qed

\vspace{0.3cm}

We next consider Case (ii), i.e.,  $a\neq 0$ and $b\not=0$. 
\begin{lemma}\label{lem:abnot0}
If $a\neq 0$ and $b\not=0$, then 
\begin{align}
V_2+V_3=&\,\frac{1}{4(q+1)}\sum_{u=0,1}\sum_{h=0}^{2q+1}G_{q^2}(\chi_{2(q+1)}^{-h}\rho^u)G_{q^2}(\chi_{2(q+1)}^{-h})\chi_{2(q+1)}^h(ab)\rho^u(a)\nonumber
\\
&\hspace{3.3cm}\times \Big(\sum_{i\in A}\chi_{2(q+1)}^{h}(\omega^i)+\sum_{i\in B}\chi_{2(q+1)}^h(\omega^{i})\rho^u(\omega)\Big),  \label{eigen34}
\end{align} 
where $\chi_{2(q+1)}$ is a  multiplicative character of order $2(q+1)$ of $\F_{q^2}$ and $\rho$ is the quadratic character of $\F_{q^2}$. 
\end{lemma}
\proof 
Let $\chi$ be a multiplicative character of order $q^2-1$ of $\F_{q^2}$. 
By \eqref{eq:ortho1}, we have 
\begin{align}
V_2=&\,\frac{1}{2(q^2-1)^2}\sum_{i\in A}\sum_{y\in C_0^{(2,q^2)}}\sum_{x\in \F_q^\ast}\sum_{j,k=0}^{q^2-2}G_{q^2}(\chi^{-j})\chi^j(axy)G_{q^2}(\chi^{-k})\chi^k(bxy^{-1}\omega^i)\nonumber\\
=&\,\frac{1}{2(q^2-1)^2}\sum_{i\in A}\sum_{j,k=0}^{q^2-2}G_{q^2}(\chi^{-j})G_{q^2}(\chi^{-k})\chi^j(a)\chi^k(b\omega^i)\chi^{j-k}(C_0^{(2,q^2)})\Big(\sum_{x\in \F_q^\ast}\chi^{j+k}(x)\Big).
\label{eigen2}
\end{align}
Since
$\chi^{j-k}(C_0^{(2,q^2)})=\frac{q^2-1}{2}$ or $0$ according as $j-k\equiv 0\,(\mod{\frac{q^2-1}{2}})$ or not, continuing from (\ref{eigen2}), we have
{\footnotesize\begin{align}
V_2=&\,\frac{1}{4(q^2-1)}\sum_{i\in A}\sum_{u=0,1}\sum_{k=0}^{q^2-2}G_{q^2}(\chi^{-k-\frac{q^2-1}{2}u})G_{q^2}(\chi^{-k})\chi^{k+\frac{q^2-1}{2}u}(a)\chi^{k}(b\omega^i)\Big(\sum_{x\in \F_q^\ast}\chi^{2k+\frac{q^2-1}{2}u}(x)\Big).\label{eigen3}
\end{align}}
Let $\chi_{2(q+1)}=\chi^{\frac{q-1}{2}}$ and $\rho=\chi^{\frac{q^2-1}{2}}$. 
Since
$\sum_{x\in \F_q^\ast}\chi^{2k+\frac{q^2-1}{2}u}(x)=q-1$ or $0$ according as $2k\equiv 0\,(\mod{q-1})$ or not, continuing from (\ref{eigen3}), we have
\[
V_2=\frac{1}{4(q+1)}\sum_{u=0,1}\sum_{h=0}^{2q+1}G_{q^2}(\chi_{2(q+1)}^{-h}\rho^u)G_{q^2}(\chi_{2(q+1)}^{-h})\chi_{2(q+1)}^h(ab)\rho^u(a)\sum_{i\in A}\chi_{2(q+1)}^{h}(\omega^i). 
\]
Similarly, we have 
\[
V_3
=\frac{1}{4(q+1)}\sum_{u=0,1}\sum_{h=0}^{2q+1}G_{q^2}(\chi_{2(q+1)}^{-h}\rho^u)G_{q^2}(\chi_{2(q+1)}^{-h})\chi_{2(q+1)}^h(ab)\rho^u(a)
\sum_{i\in B}\chi_{2(q+1)}^h(\omega^{i})\rho^u(\omega).
\]
This completes the proof of the lemma. \qed 
\vspace{0.3cm}

Let $W_0$ (resp. $W_1$) be the contribution for $u=0$ (resp. $u=1$) in the summations of (\ref{eigen34}); then $V_2+V_3=W_0+W_1$.  \begin{lemma}\label{lem:ss1}
Let $r=ab\not=0$. Then, 
\[
W_0=
\left\{
\begin{array}{ll}
\frac{-q^2+1}{4}, & \mbox{ if $r\in \omega^c R_1$\, or\,  $r\in \omega^c R_2$,}\\
\frac{q^2+1}{4}\mbox{ or } \frac{-3q^2+1}{4}, &\mbox{ otherwise, }
 \end{array}
\right.
\] 
depending on  whether $q\equiv 3$ or $1\,(\mod{4})$.
\end{lemma}
\proof 
By the definition of $W_0$, we have 
\begin{align*}
W_0=\frac{1}{4(q+1)}\sum_{h=0}^{2q+1}G_{q^2}(\chi_{2(q+1)}^{-h})^2\chi_{2(q+1)}^h(r)\Big(\sum_{i\in A\cup B}\chi_{2(q+1)}^{h}(\omega^i)\Big).  
\end{align*}
Since $A\cup B=X_{1,c}\cup X_{3,c}\cup X_{3,c}$ as a multiset, by the property~(P7), we have 
\[
\psi_{\F_{q^2}}(\omega^a\bigcup_{i\in A\cup B}C_i^{(2(q+1),q^2)})=\left\{
\begin{array}{ll}
\frac{q G_q(\eta)-1}{2}(=:c_1), & \mbox{ if $a\in Y_{1,c}(=:Z_1)$,}\\
\frac{- q G_q(\eta)-1}{2}(=:c_2), & \mbox{ if $a\in Y_{2,c}(=:Z_2)$,}\\
\frac{-1+(-1)^{\frac{q-1}{2}}q}{2}(=:c_3), & \mbox{ if $a\in Y_{3,c}\cup Y_{4,c}(=:Z_3)$,}\\
\frac{-1-(-1)^{\frac{q-1}{2}}q}{2}(=:c_4), & \mbox{ if $a\in Y_{5,c}(=:Z_4)$. }
 \end{array}
\right.
\]
Then, by \eqref{eq:ortho2}, we have 
\begin{align*}
G_{q^2}(\chi_{2(q+1)}^{-h})\sum_{i\in A\cup B}\chi_{2(q+1)}^{h}(\omega^{i})
=&\,
\sum_{a=0}^{2q+1}\psi_{\F_{q^2}}(\omega^a \bigcup_{i\in A\cup B}C_i^{(2(q+1),q^2)})
\chi_{2(q+1)}^{-h}(\omega^a)\\
=&\,\sum_{i=1}^4c_{i}\sum_{a\in Z_i}\chi_{2(q+1)}^{-h}(\omega^a). 
\end{align*}
Then, by \eqref{eq:ortho1},  
we have 
\begin{align*}
W_0=&\,\sum_{i=1}^4\frac{c_i}{4(q+1)}\sum_{a\in Z_i}\sum_{h=0}^{2q+1}G_{q^2}(\chi_{2(q+1)}^{-h})\chi_{2(q+1)}^h(r\omega^{-a})\\
=&\,\sum_{i=1}^4\frac{c_i}{2}
\sum_{a\in -Z_{i}}\psi_{\F_{q^2}}(rC_a^{(2(q+1),q^2)})
.
\end{align*}
Since $-Y_{i,c}\equiv Y_{i,c}-c\,(\mod{2(q+1)})$, $i=1,2$, from the property~(P8),  we have by the property~(P10) that for $i=1,2$   
\[
\sum_{a\in -Z_{i}}\psi_{\F_{q^2}}(rC_a^{(2(q+1),q^2)})
=
\left\{
\begin{array}{ll}
\frac{-2+q+(-1)^{i-1} G_q(\eta)}{2}, & \mbox{ if $r\in \omega^cR_1$,}\\
\frac{-2+q-(-1)^{i-1} G_q(\eta)}{2}, &\mbox{ if  $r\in \omega^cR_2$,}\\
-1+(-1)^{i-1} G_q(\eta), & \mbox{ if  $r\in \omega^cR_3$,}\\
-1-(-1)^{i-1} G_q(\eta),& \mbox{ if  $r\in \omega^cR_4$,}\\
-1,& \mbox{ if  $r\in\omega^cR_5$.} 
 \end{array}
\right.
\]
Furthermore, since $-(Y_{3,c}\cup Y_{4,c})\equiv Y_{5,c}-c\,(\mod{2(q+1)})$ by the property~(P9), we have 
\begin{align*}
\sum_{a\in -Z_3}\psi_{\F_{q^2}}(rC_a^{(2(q+1),q^2)})=&\, 
\sum_{a\in Y_{5,c}-c}\psi_{\F_{q^2}}(rC_a^{(2(q+1),q^2)})
\\
=&\, 
\left\{
\begin{array}{ll}
\frac{1-q}{2}, & \mbox{ if $r\in \omega^c(R_1\cup R_2)$,}\\
\frac{1+(-1)^{\frac{q-1}{2}}q}{2}, &\mbox{ if $r\in \omega^c R_5$,}\\
\frac{1-(-1)^{\frac{q-1}{2}}q}{2}, &\mbox{ if $r\in \omega^c(R_3\cup R_4)$.}
 \end{array}
\right. 
\end{align*}
Similarly,  
we have 
\begin{align*}
\sum_{a\in -Z_4}\psi_{\F_{q^2}}(rC_a^{(2(q+1),q^2)})=&\, 
\sum_{a\in (Y_{3}\cup Y_{4})-c}\psi_{\F_{q^2}}(rC_a^{(2(q+1),q^2)})\\
=&\,
\left\{
\begin{array}{ll}
\frac{-q+1}{2}, & \mbox{ if $r\in \omega^c(R_1\cup R_2)$,}\\
\frac{1-(-1)^{\frac{q-1}{2}}q}{2}, & \mbox{ if $r\in \omega^c R_5$,}\\
\frac{1+(-1)^{\frac{q-1}{2}}q}{2}, & \mbox{ if $r\in \omega^c(R_3\cup R_4)$.}
 \end{array}
\right.
\end{align*}
Summing up, we have 
\[
W_0=
\left\{
\begin{array}{ll}
\frac{-q^2+1}{4}, & \mbox{ if $r\in \omega^c R_1$\, or\,  $r\in \omega^c R_2$,}\\
\frac{q^2+1}{4}, &\mbox{ if $r\in \omega^c (R_2\cup R_4\cup R_5)$\, or\, $r\in \omega^c (R_1\cup R_3\cup R_5)$,}\\
\frac{-3q^2+1}{4}, & \mbox{ if $r\in \omega^c R_3$\, 
or\, $r\in \omega^c R_4$, }
 \end{array}
\right.
\] 
according as $q\equiv 3$ or $1\,(\mod{4})$. This completes the proof. 
\qed
\vspace{0.3cm}

We next evaluate $W_1$ below.   
\begin{lemma}\label{lem:s1}
Let $r=ab\not=0$. Then, 
\begin{align*}
W_1=&\,
-\frac{(-1)^{\frac{q-1}{2}}\rho(a)q}{4}\Big(|A|-|B|+\rho(r)(|A_e|+|B_o|-|A_o|-|B_e|)\Big)\\
&\quad +
\frac{(-1)^{\frac{q-1}{2}}\rho(a)q^2}{2}\cdot
\left\{
\begin{array}{ll}
1,& \mbox{ if  $r \in \bigcup_{i\in -(A\setminus B)+q+1}C_{i}^{(2(q+1),q^2)}$,} \\
-1,& \mbox{ if $r \in \bigcup_{i\in -(B\setminus A)+q+1}C_{i}^{(2(q+1),q^2)}$,} \\
0,& \mbox{ otherwise. }
 \end{array}
\right.
\end{align*}
\end{lemma}
\proof 
By the definition of $W_1$, we have 
\begin{equation}\label{eq:w1}
W_1=\frac{\rho(a)}{4(q+1)}\sum_{h=0}^{2q+1}G_{q^2}(\chi_{2(q+1)}^{-h}\rho)G_{q^2}(\chi_{2(q+1)}^{-h})\chi_{2(q+1)}^h(r)\Big(\sum_{i\in A}\chi_{2(q+1)}^{h}(\omega^i)-\sum_{i\in B}\chi_{2(q+1)}^h(\omega^{i})\Big).
\end{equation}
By applying the Davenport-Hasse product  formula (Theorem~\ref{thm:Stickel2}) with $\chi=\chi_{2(q+1)}^{-h}$, $\chi'=\rho$, and $\ell=2$ we have 
\[
G_{q^2}(\chi_{2(q+1)}^{-h})G_{q^2}(\chi_{2(q+1)}^{-h}\rho)=G_{q^2}(\rho)G_{q^2}(\chi_{q+1}^{-h}),  
\]
where $\chi_{q+1}=\chi_{2(q+1)}^{2}$ has order $q+1$. Then, \eqref{eq:w1} is rewritten as 
\begin{equation}\label{eq:s0}
W_1=\frac{\rho(a)}{4(q+1)}G_{q^2}(\rho)\sum_{h=0}^{2q+1}G_{q^2}(\chi_{q+1}^{-h})\chi_{2(q+1)}^h(r)\Big(\sum_{i\in A}\chi_{2(q+1)}^{h}(\omega^i)-\sum_{i\in B}\chi_{2(q+1)}^h(\omega^{i})\Big). 
\end{equation}
We will compute $W_1$ by dividing it into three parts. Let $W_{1,1},W_{1,2},W_{1,3}$ denote the contributions in the sum on the right hand side of (\ref{eq:s0}) when $h=0, q+1$;  other even $h$; and odd $h$, respectively. Then $W_1=W_{1,1}+W_{1,2}+W_{1,3}$. For $W_{1,1}$, we have  
\[
W_{1,1}=-\frac{\rho(a)}{4(q+1)}G_{q^2}(\rho)\Big(|A|-|B|+\rho(r)(|A_e|+|B_o|-|A_o|-|B_e|)\Big).
\]
Next, by Theorem~\ref{thm:semi}, we have 
\begin{equation}\label{eq:w12}
W_{1,2}=\frac{\rho(a)q}{4(q+1)}G_{q^2}(\rho)\sum_{\ell=0;\ell\not=0,\frac{q+1}{2}}^{q}\chi_{q+1}^\ell(r)\Big(\sum_{i\in A}\chi_{q+1}^{\ell}(\omega^i)-\sum_{i\in B}\chi_{q+1}^\ell(\omega^{i})\Big).
\end{equation}
By the property~(P2), 
\[
\{x\,(\mod{q+1})\,|\,x\in A\setminus B\}\cap \{x\,(\mod{q+1})\,|\,x\in B\setminus A\}=\emptyset. 
\]
Hence, 
continuing from \eqref{eq:w12}, we have
\begin{align*} 
W_{1,2}=& \,-\frac{\rho(a)q}{4(q+1)}G_{q^2}(\rho)\Big(|A|-|B|+\rho(r)(|A_e|+|B_o|-|A_o|-|B_e|)\Big)\\
&\quad +\frac{\rho(a)q}{4}G_{q^2}(\rho)
\cdot 
\left\{
\begin{array}{ll}
1,& \mbox{ if $r \in \bigcup_{i\in -(A\setminus B)}C_{i}^{(q+1,q^2)}$,} \\
-1,& \mbox{ if $r \in \bigcup_{i\in -(B\setminus A)}C_{i}^{(q+1,q^2)}$,}\\
0,& \mbox{ otherwise. } \\
 \end{array}
\right.
\end{align*}
Finally, by Theorem~\ref{thm:semi} again, we have 
\begin{align*}
W_{1,3}=&\,-\frac{\rho(a)q}{4(q+1)}G_{q^2}(\rho)\sum_{h:\, odd}\chi_{2(q+1)}^h(r)\Big(\sum_{i\in A}\chi_{2(q+1)}^{h}(\omega^i)-\sum_{i\in B}\chi_{2(q+1)}^h(\omega^{i})\Big)\\
=&\,-\frac{\rho(a)q}{4(q+1)}G_{q^2}(\rho)\sum_{h=0}^{2q+1}\chi_{2(q+1)}^h(r)\Big(\sum_{i\in A}\chi_{2(q+1)}^{h}(\omega^i)-\sum_{i\in B}\chi_{2(q+1)}^h(\omega^{i})\Big)\\
&\, \quad +\frac{\rho(a)q}{4(q+1)}G_{q^2}(\rho)\sum_{\ell=0}^{q}\chi_{q+1}^\ell(r)\Big(\sum_{i\in A}\chi_{q+1}^{\ell}(\omega^i)-\sum_{i\in B}\chi_{q+1}^\ell(\omega^{i})\Big)\\
=&\,-\frac{\rho(a)q}{2}G_{q^2}(\rho)\cdot 
\left\{
\begin{array}{ll}
1,& \mbox{ if $r \in \bigcup_{i\in -(A\setminus B)}C_{i}^{(2(q+1),q^2)}$,} \\
-1,& \mbox{ if $r \in \bigcup_{i\in -(B\setminus A)}C_{i}^{(2(q+1),q^2)}$,}\\
0,& \mbox{ otherwise, } 
 \end{array}
\right.\\
&\quad +\frac{\rho(a)q}{4}G_{q^2}(\rho)
\cdot \left\{
\begin{array}{ll}
1,& \mbox{ if $r \in \bigcup_{i\in -(A\setminus B)}C_{i}^{(q+1,q^2)}$,} \\
-1,& \mbox{ if $r \in \bigcup_{i\in -(B\setminus A)}C_{i}^{(q+1,q^2)}$,}\\
0,& \mbox{ otherwise. } 
 \end{array}
\right.
\end{align*}
Summing up, we have 
\begin{align*}
W_1=&\,W_{1,2}+W_{1,2}+W_{1,3}\\
=&\,
-\frac{\rho(a)}{4}G_{q^2}(\rho)\Big(|A|-|B|+\rho(r)(|A_e|+|B_o|-|A_o|-|B_e|)\Big)\\
&\quad +
\frac{\rho(a)q}{2}G_{q^2}(\rho)\cdot
\left\{
\begin{array}{ll}
1,& \mbox{ if  $r \in \bigcup_{i\in -(A\setminus B)+q+1}C_{i}^{(2(q+1),q^2)}$,} \\
-1,& \mbox{ if $r \in \bigcup_{i\in -(B\setminus A)+q+1}C_{i}^{(2(q+1),q^2)}$,} \\
0,& \mbox{ otherwise. }
 \end{array}
\right.
\end{align*}
The statement now follows from $G_{q^2}(\rho)=-(-1)^{\frac{q-1}{2}}q$.  
\qed

\begin{remark}
By Lemmas~\ref{lem:ss1} and ~\ref{lem:s1}, we have 
\begin{align}
V_2+V_3=&\, W_0+W_1\nonumber\\
=&\, 
\frac{(-1)^{\frac{q-1}{2}}\rho(a)q}{4}\Big(|A|-|B|+\rho(r)(|A_e|+|B_o|-|A_o|-|B_e|)\Big)\nonumber\\
&\, \quad -\frac{(-1)^{\frac{q-1}{2}}\rho(a)q^2}{2}\cdot \left\{
\begin{array}{ll}
1,& \mbox{ if  $r \in \bigcup_{i\in -(A\setminus B)+(q+1)}C_i^{(2(q+1),q^2)}$,} \\
-1,& \mbox{ if $r \in \bigcup_{i\in -(B\setminus A)+(q+1)}C_i^{(2(q+1),q^2)}$,} \\
0,& \mbox{ otherwise,}
 \end{array}
\right.\nonumber\\
&\, \, \quad +\left\{
\begin{array}{ll}
\frac{-q^2+1}{4},& \mbox{ if $r\in \omega^c R_1$ or $r\in \omega^c R_2$,} \\
\frac{-3q^2+1}{4} \mbox{\, or \,} \frac{q^2+1}{4},& \mbox{ otherwise,}
 \end{array}
\right.
\label{eq:chara}
\end{align}
according as $q\equiv 3$ or $1\,(\mod{4})$. 
By the property~(P4), $X_{1,c}+c\equiv -((A\setminus B)\cup (B\setminus A))+(q+1)\,(\mod{2(q+1)})$ or  $X_{2,c}+c\equiv-((A\setminus B)\cup (B\setminus A))+(q+1)\,(\mod{2(q+1)})$ depending on whether $q\equiv 3$ or $1\,(\mod{4})$. Hence, continuing from \eqref{eq:chara}, we have  
\begin{align}
&V_2+V_3-\frac{(-1)^{\frac{q-1}{2}}\rho(a)q}{4}\Big(|A|-|B|+\rho(r)(|A_e|+|B_o|-|A_o|-|B_e|)\Big) \nonumber\\
=&\, \frac{q^2+1}{4}\, \, \mbox{or }\, \, \frac{-3q^2+1}{4}.\label{eq:v2v3last}
\end{align}
\end{remark}
\vspace{0.3cm}

We are now ready to prove our main theorem. 

\vspace{0.1cm}
{\bf Proof of Theorem~\ref{thm:mainWX1}:\,} 
In the case where exactly one of $a, b$ is zero, the statement follows from 
Lemmas~\ref{lemma:ab0} and \ref{lemma:ab02}. We treat the case where $a\neq 0$ and $b\not=0$. 

(1)  By \eqref{eq:Gaussquad}, $V_0=\frac{-1+\rho(b)q}{2}$. 
Furthermore, 
by $|A|=|B|=\frac{q+1}{2}$ and $|A_e|+|B_o|=|A_o|+|B_e|-2(-1)^{\frac{q-1}{2}}$, 
we have 
\[
\frac{(-1)^{\frac{q-1}{2}}\rho(a)q}{4}\Big(|A|-|B|+\rho(r)(|A_e|+|B_o|-|A_o|-|B_e|)\Big)=-\frac{\rho(b)q}{2}. 
\]
Hence, by \eqref{eq:v2v3last},  it follows that 
$\psi_{a,b}(E_0)=
\frac{q^2-1}{4}\, \, \mbox{or }\, \, \frac{-3q^2-1}{4}$.

(2) By \eqref{eq:Gaussquad}, $V_1=\frac{-1-(-1)^{\frac{q-1}{2}}\rho(a)q}{2}$. 
Furthermore, 
by $|A|=\frac{q+3}{2}$, 
$|B|=\frac{q-1}{2}$, and $|A_e|+|B_o|=|A_o|+|B_e|$, 
we have 
\[
\frac{(-1)^{\frac{q-1}{2}}\rho(a)q}{4}\Big(|A|-|B|+\rho(r)(|A_e|+|B_o|-|A_o|-|B_e|)\Big)=\frac{(-1)^{\frac{q-1}{2}}\rho(a)q}{2}. 
\]
Hence, by \eqref{eq:v2v3last},  it follows that 
$\psi_{a,b}(E_1)=
\frac{q^2-1}{4}\, \, \mbox{or }\, \, \frac{-3q^2-1}{4}$.
\qed

\begin{corollary}
Let $A=\{\beta\}\cup X_{3,c}$, $B=\{\alpha\}\cup X_{3,c}$, 
$A'=X_{1,c+q+1}\cup X_{3,c+q+1}$, and $B'=X_{3,c+q+1}$, where $\alpha,\beta$ are defined as in the property~(P6).  
Then, 
the sets
\begin{align*}
C_0=&\{(0,y)\,|\,y\in C_\tau^{(2,q^2)}\} \cup \{(xy,xy^{-1}\omega^i)\,|\,x\in \F_q^\ast,y\in C_0^{(2,q^2)},i\in A\}\\
&\quad \cup  
\{(xy,xy^{-1}\omega^i)\,|\,x\in \F_q^\ast,y\in C_1^{(2,q^2)},i\in B\},\\
C_1=&\{(y,0)\,|\,y\in C_0^{(2,q^2)}\}\cup \{(xy,xy^{-1}\omega^i)\,|\,x\in \F_q^\ast,y\in C_0^{(2,q^2)},i\in A'\}\\
&\quad \cup  \{(xy,xy^{-1}\omega^i)\,|\,x\in \F_q^\ast,y\in C_1^{(2,q^2)},i\in B'\},\\
C_2=&\{(0,y)\,|\,y\in C_{\tau+1}^{(2,q^2)}\} \cup \{(xy\omega,xy^{-1}\omega^{i+q})\,|\,x\in \F_q^\ast,y\in C_0^{(2,q^2)},i\in A\}\\
&\quad \cup  
\{(xy\omega,xy^{-1}\omega^{i+q})\,|\,x\in \F_q^\ast,y\in C_1^{(2,q^2)},i\in B\},\\
C_3=&\{(y,0)\,|\,y\in C_1^{(2,q^2)}\}\cup \{(xy\omega,xy^{-1}\omega^{i+q})\,|\,x\in \F_q^\ast,y\in C_0^{(2,q^2)},i\in A'\}\\
&\quad \cup  \{(xy\omega,xy^{-1}\omega^{i+q})\,|\,x\in \F_q^\ast,y\in C_1^{(2,q^2)},i\in B'\}
\end{align*} 
are of type Q.  Furthermore, these sets satisfy the assumptions of Remark~\ref{rem:HDiff} with respect to the spread ${\mathcal K}$ consisting of the following $2$-dimensional subspaces: 
\[
K_y=\{ (x,yx^q) \,|\,x\in \F_{q^2}\}, y\in \F_{q^2}, \, \mbox{\, and \, }\, K_\infty=\{(0,x)\,|\,x\in \F_{q^2}\}. 
\]
\end{corollary}
\proof 
By the property (P6), $|A_e|+|B_o|=|A_o|+|B_e|-2(-1)^{\frac{q-1}{2}}$ and 
$|A_e'|+|B_o'|=|A_o'|+|B_e'|$. 
Hence, by 
Theorem~\ref{thm:mainWX1}, $C_0$ and $C_1$ are type Q sets. Since $C_2=\{(\omega x,\omega^q y)\,|\, (x,y)\in C_0\}$ and $C_3=\{(\omega x,\omega^q y)\,|\, (x,y)\in C_1\}$, the sets $C_2$ and $C_3$ are also of type Q. Furthermore, $\bigcup_{i=1}^3C_i=(\F_{q^2}\times \F_{q^2})\setminus \{(0,0)\}$ since $A\cup A'\cup (B+q+1) \cup (B'+q+1)\equiv \{0,1,\ldots,2q+1\}\,(\mod{2(q+1)})$ by the properties~(P1),(P4),(P5) and (P6).  Therefore,  
 $C_i$, $i=0,1,2,3$, satisfy 
the assumptions of Remark~\ref{rem:HDiff}  as 
$
C_0\cup C_2\cup \{(0,0)\}=\big(\bigcup_{y\in H_0}K_y\big)\cup K_\infty$ and 
$C_1\cup C_3\cup \{(0,0)\}=\bigcup_{y\in H_1}K_y$, where $H_0= \bigcup_{i\in -(A\cup (B+q+1))}C_i^{(2(q+1),q^2)}$ and $H_1=\bigcup_{i\in -(A'\cup (B'+q+1))}C_i^{(2(q+1),q^2)}$.  
\qed

\section*{Appendix}\label{sec:newass}
In this appendix, we prove that the sets $X_{i,c}$ and $Y_{i,c}$, $i=1,2,3,4,5$, have the properties (P1)--(P10).  

By the definition of $X_{1,c}$, we have 
\begin{align*}
X_{1,c}=&\,(\{\tfrac{q+1}{2},\tfrac{3(q+1)}{2}\}\cap \{i \, (\mod{2(q+1)})\,|\,\Tr_{q^2/q}(\omega^{i+c})\in C_0^{(2,q)}\})\\
&\quad \cup (\{\tfrac{q+1}{2}-c,\tfrac{3(q+1)}{2}-c\}\cap 
\{i \, (\mod{2(q+1)})\,|\,\Tr_{q^2/q}(\omega^{i})\in C_0^{(2,q)}\}). 
\end{align*}
Hence, there are $\epsilon,\delta\in \{-1,1\}$ such that 
$X_{1,c}=\{\frac{q+1}{2}\epsilon,\frac{q+1}{2}\delta-c\}$. In particular, 
we have 
\begin{equation}\label{eq:traceapp1}
\Tr_{q^2/q}(\omega^{c+\frac{q+1}{2}\epsilon})\in C_0^{(2,q)}
\mbox{\, \, and\, \, }
\Tr_{q^2/q}(\omega^{\frac{q+1}{2}\delta-c})\in C_0^{(2,q)}. 
\end{equation}
\begin{lemma}\label{lem:x1c}
We have $X_{1,c}=\{\tfrac{q+1}{2}\epsilon,\tfrac{q+1}{2}\delta-c\}$ 
for some $(\epsilon,\delta)\in \{(1,1),(-1,-1)\} $ or $\{(-1,1),(1,-1)\} $ according as $q\equiv 3\,(\mod{4})$ or $q\equiv 1\,(\mod{4})$.  
\end{lemma}
\proof 
By \eqref{eq:traceapp1}, we have 
\begin{equation}\label{eq:traceapp0}
\omega^{\frac{q+1}{2}\epsilon}(\omega^{c}-\omega^{cq})\in C_0^{(2,q)}\mbox{\, \, and\, \, }
\omega^{\frac{q+1}{2}\delta}(\omega^{-c}-\omega^{-cq})\in C_0^{(2,q)}. 
\end{equation}
Putting $\omega^d=1-\omega^{c(q-1)}$, the conditions in \eqref{eq:traceapp0} are rewritten as 
\[
\omega^{\frac{q+1}{2}\epsilon+c+d}= \omega^{2(q+1)k}\mbox{\, \, and\, \, }
\omega^{\frac{q+1}{2}\delta-c+dq}=\omega^{2(q+1)\ell} 
\]
for some $k,\ell\in \Z$. Here, $d$ is odd if $q\equiv 3\,(\mod{4})$, and  $d$ is even if $q\equiv 1\,(\mod{4})$. By multiplying these 
equations, we have $\omega^{\frac{q+1}{2}(\epsilon+\delta)+d(q+1)}=\omega^{2(q+1)(k+\ell)}$. Then, the statement immediately follows. 
\qed
\begin{remark}\label{rem:fives}
For  $X_{i,c}$, $i=1,2,3,4,5$, we observe the following facts:  
\begin{itemize}
\item[(1)]
Since $I_2\equiv I_3+(q+1)\,(\mod{2(q+1)})$, we have 
$X_{1,c}\equiv X_{2,c}+(q+1)\,(\mod{2(q+1)})$, 
$X_{3,c}\equiv X_{4,c}+(q+1)\,(\mod{2(q+1)})$, and $X_{5,c}\equiv X_{5,c}+(q+1)\,(\mod{2(q+1)})$. Hence, the property~(P1) follows. 
\item[(2)]
Since $I_2$ forms a $(q+1,2,q,\frac{q-1}{2})$ relative difference set  (cf.~\cite{ADJP}), 
we have 
$|X_{3,c}|=\frac{q-1}{2}$. Then, the property~(P2) follows. 
\item[(3)]
Since $X_{3,c+q+1}=I_2\cap J_3$ and  $X_{4,c+q+1}=I_3\cap J_2$, we have 
$X_{3,c+q+1}\cup X_{4,c+q+1}=X_{5,c}$. Then, the property~(P3) follows. 
\item[(4)] 
The property~(P4) directly follows from Lemma~\ref{lem:x1c}. 
\item[(5)] By Lemma~\ref{lem:x1c}, 
$X_{1,c+q+1}=\{\frac{q+1}{2}\epsilon',\frac{q+1}{2}\delta'-c+q+1\}$
for some $(\epsilon',\delta')\in \{(1,1),(-1,-1)\} $ or $\{(-1,1),(1,-1)\} $ according to whether $q\equiv 3\,(\mod{4})$ or $q\equiv 1\,(\mod{4})$. Then, it is direct to see that $|X_{1,c}\cap X_{1,c+q+1}|=1$ 
and $(X_{1,c}\setminus X_{1,c+q+1})\equiv (X_{1,c+q+1}\setminus X_{1,c})+q+1\,(\mod{2(q+1)})$ in all cases. More precisely, $X_{1,c+q+1}=\{\frac{q+1}{2}\epsilon+q+1,\frac{q+1}{2}\delta-c\}$ since $\frac{q+1}{2}\delta-c\in J_1\cap I_2$. Hence, $X_{1,c}\setminus X_{1,c+q+1}=\{\frac{q+1}{2}\epsilon\}$ and $X_{1,c}\cap  X_{1,c+q+1}=\{\frac{q+1}{2}\delta-c\}$. Thus, the properties~(P5) and (P6) follow. 
\end{itemize} 
\end{remark}
Next, we show that the $X_{i,c}$'s have property~(P7). 
\begin{proposition}\label{prop:asso}
Let $R_{i}$, $i=1,2,3,4,5$, be defined as  in Subsection~\ref{subsec:set}. Then, $R_i$, $i=1,2,3,4,5$, take the character values listed in Table~\ref{tab_1}. In particular, 
$Y_{i,c}$'s in Table~\ref{tab_1} are determined as follows: 
\begin{align*}
Y_{1,c}=&\{0,c\}, \, \, \, Y_{2,c}=\{q+1,c+q+1\},\\
Y_{3,c}=&\{i+c-\tfrac{q+1}{2}\delta\,|\,\Tr_{q^2/q}(\omega^i)\in C_0^{(2,q)}\}\cap \{i-\tfrac{q+1}{2}\epsilon\,|\,\Tr_{q^2/q}(\omega^i)\in C_0^{(2,q)}\}, \\
Y_{4,c}=&\{i+c-\tfrac{q+1}{2}\delta\,|\,\Tr_{q^2/q}(\omega^i)\in C_1^{(2,q)}\}\cap \{i-\tfrac{q+1}{2}\epsilon\,|\,\Tr_{q^2/q}(\omega^i)\in C_1^{(2,q)}\},\\
Y_{5,c}=& (\{i+c-\tfrac{q+1}{2}\delta\,|\,\Tr_{q^2/q}(\omega^i)\in C_0^{(2,q)}\}\cap \{i-\tfrac{q+1}{2}\epsilon\,|\,\Tr_{q^2/q}(\omega^i)\in C_1^{(2,q)}\})\\
&\, \cup  (\{i+c-\tfrac{q+1}{2}\delta\,|\,\Tr_{q^2/q}(\omega^i)\in C_1^{(2,q)}\}\cap \{i-\tfrac{q+1}{2}\epsilon\,|\,\Tr_{q^2/q}(\omega^i)\in C_0^{(2,q)}\}). 
\end{align*}
\end{proposition}
\proof
The character values $\psi_{\F_{q^2}}(\omega^aR_1)$, $a=0,1,\ldots,2q+1$, are evaluated as follows: 
\begin{align*}
\psi_{\F_{q^2}}(\omega^aR_1)=&\,\psi_{\F_{q^2}}(\omega^{a+\frac{q+1}{2}\delta-c}C_0^{(2(q+1),q^2)})+\psi_{\F_{q^2}}(\omega^{a+\frac{q+1}{2}\epsilon}C_0^{(2(q+1),q^2)})\\
=&\,
\psi_{\F_q}(\Tr_{q^2/q}(\omega^{a+\frac{q+1}{2}\delta-c})C_0^{(2,q)})+
\psi_{\F_q}(\Tr_{q^2/q}(\omega^{a+\frac{q+1}{2}\epsilon})C_0^{(2,q)})
\\
=&\,
\left\{
\begin{array}{ll}
\frac{q-1}{2}, & \mbox{ if $a\in I_1-\frac{q+1}{2}\delta+c$,}\\
\frac{-1+G_q(\eta)}{2}, &\mbox{ if $a\in I_2-\frac{q+1}{2}\delta+c$,}\\
\frac{-1-G_q(\eta)}{2}, & \mbox{ if $a\in I_3-\frac{q+1}{2}\delta+c$,}
 \end{array}
\right.
+
\left\{
\begin{array}{ll}
\frac{q-1}{2}, & \mbox{ if $a\in I_1-\frac{q+1}{2}\epsilon$,}\\
\frac{-1+G_q(\eta)}{2}, &\mbox{ if $a\in I_2-\frac{q+1}{2}\epsilon$,}\\
\frac{-1-G_q(\eta)}{2}, & \mbox{ if $a\in I_3-\frac{q+1}{2}\epsilon$, }
 \end{array}
\right.\\
=&\,
\left\{
\begin{array}{ll}
\frac{-2+q+G_q(\eta)}{2}, & \mbox{ if $a\in Y_{1,c}'$,}\\
\frac{-2+q-G_q(\eta)}{2}, &\mbox{ if $a\in Y_{2,c}'$,}\\
-1+G_q(\eta), & \mbox{ if $a\in Y_{3,c}$,}\\
-1-G_q(\eta),& \mbox{ if $a\in Y_{4,c}$,}\\
-1,& \mbox{ if $a\in Y_{5,c}$, } 
 \end{array}
\right.
\end{align*}
where
\begin{align*}
Y_{1,c}'=&\,((I_1-\tfrac{q+1}{2}\delta +c)\cap (I_2-\tfrac{q+1}{2}
\epsilon)) \cup 
((I_2-\tfrac{q+1}{2}\delta +c)\cap (I_1-\tfrac{q+1}{2}\epsilon)),\\
Y_{2,c}'=&\,((I_1-\tfrac{q+1}{2}\delta +c)\cap (I_3-\tfrac{q+1}{2}\epsilon))\cup 
((I_3-\tfrac{q+1}{2}\delta +c)\cap (I_1-\tfrac{q+1}{2}\epsilon)). 
\end{align*}
By (\ref{eq:traceapp1}), it is direct to see that 
\begin{align*}
&(I_1-\tfrac{q+1}{2}\delta +c)\cap (I_2-\tfrac{q+1}{2}
\epsilon)=\{c\},\quad 
(I_2-\tfrac{q+1}{2}\delta +c)\cap (I_1-\tfrac{q+1}{2}\epsilon)=\{0\}, \\
&(I_1-\tfrac{q+1}{2}\delta +c)\cap (I_3-\tfrac{q+1}{2}
\epsilon)=\{c+q+1\},\quad 
(I_3-\tfrac{q+1}{2}\delta +c)\cap (I_1-\tfrac{q+1}{2}\epsilon)=\{q+1\}. 
\end{align*}
Hence, we have $Y_{1,c}'=Y_{1,c}$ and $Y_{2,c}'=Y_{2,c}$. 

The character values of $R_2$ is determined as $\psi_{\F_{q^2}}(\omega^a R_2)=\psi_{\F_{q^2}}(\omega^{a+q+1} R_1)$. 

We next evaluate $\psi_{\F_{q^2}}(\omega^a R_3)$, $a=0,1,\ldots,2q+1$. 
By Remark~\ref{rem:secChen}~(i), the indicator function of 
$\{x\,|\,\Tr_{q^2/q}(x)\in C_0^{(2,q)}\}$ is given by 
\[
f(x)=\frac{1}{q}\sum_{s\in \F_q}\sum_{y\in C_0^{(2,q)}}
\psi_{\F_{q^2}}(sx)\psi_{\F_{q}}(-sy). 
\]
Then, 
\begin{align*}
\psi_{\F_{q^2}}(\omega^a R_3)=&\,\sum_{x\in \F_{q^2}}\psi_{\F_{q^2}}(\omega^a x)f(x)f(x\omega^c)\\
=&\,\frac{1}{q^2}\sum_{x\in \F_{q^2}}\sum_{s,t\in \F_{q}}\sum_{y,z\in C_0^{(2,q)}}
\psi_{\F_{q^2}}(x(\omega^a +s+t\omega^c))\psi_{\F_{q}}(-sy)\psi_{\F_{q}}(-tz)\\
=&\,\sum_{s,t\in \F_{q}:\omega^a=s+t\omega^{c}}\sum_{y,z\in C_0^{(2,q)}}
\psi_{\F_{q}}(sy)\psi_{\F_{q}}(tz).
\end{align*}
We treat the case where $a\in Y_{1,c}\cup Y_{2,c}=\{0,c,q+1,c+q+1\}$. 
If $a=c$, then $s=0$ and $t\in C_0^{(2,q)}$, and hence $\psi_{\F_{q^2}}(\omega^aR_3)=\frac{(q-1)(-1+G_q(\eta))}{4}$. 
If  $a=c+q+1$, then $s=0$ and $t\in C_1^{(2,q)}$, and hence $\psi_{\F_{q^2}}(\omega^aR_3)=\frac{(q-1)(-1-G_q(\eta))}{4}$. 
If $a=0$, then $t=0$ and $s\in C_0^{(2,q)}$, and hence $\psi_{\F_{q^2}}(\omega^aR_3)=\frac{(q-1)(-1+G_q(\eta))}{4}$. 
If $a=q+1$, then $t=0$ and $s\in C_1^{(2,q)}$, and hence $\psi_{\F_{q^2}}(\omega^aR_3)=\frac{(q-1)(-1-G_q(\eta))}{4}$. 
Next, we treat the case where $s,t\not=0$. Define 
\begin{align*}
G_3=&\{a\,(\mod{2(q+1)})\,|\,\omega^a=s+t\omega^c,s,t\in C_0^{(2,q)}\},\\
G_4=&\{a\,(\mod{2(q+1)})\,|\,\omega^a=s+t\omega^c,s,t\in C_1^{(2,q)}\},\\
G_5=&\{a\,(\mod{2(q+1)})\,|\,\omega^a=s+t\omega^c,(s,t)\in C_0^{(2,q)}\times C_1^{(2,q)}\mbox{ or }  C_1^{(2,q)}\times C_0^{(2,q)}\}. 
\end{align*}
Then, we have 
\[
\psi_{\F_{q^2}}(\omega^aR_3)=
\left\{
\begin{array}{ll}
\frac{(1-G_q(\eta))^2}{4}, & \mbox{ if $a\in G_3$,}\\
\frac{(1+G_q(\eta))^2}{4}, &\mbox{ if $a\in G_4$, }\\
\frac{1-(-1)^{\frac{q-1}{2}}q}{4}, & \mbox{ if $a\in G_5$.}
 \end{array}
\right.
\]
We need to show that $G_i=Y_{i,c}$, $i=3,4,5$. Let $a\in G_3$. Then, there 
are some $s,t\in C_0^{(2,q)}$ such that $\omega^a=s+t\omega^c$.  
Taking trace of both sides of  $\omega^{a+\frac{q+1}{2}\epsilon}=s\omega^{\frac{q+1}{2}\epsilon}+t\omega^{c+\frac{q+1}{2}\epsilon}$, we have 
$\Tr_{q^2/q}(\omega^{a+\frac{q+1}{2}\epsilon})=
s\Tr_{q^2/q}(\omega^{\frac{q+1}{2}\epsilon})+t\Tr_{q^2/q}(\omega^{\frac{q+1}{2}\epsilon+c})$. Since $\Tr_{q^2/q}(\omega^{\frac{q+1}{2}\epsilon})=0$ and $\Tr_{q^2/q}(\omega^{\frac{q+1}{2}\epsilon+c})\in C_0^{(2,q)}$, we obtain  
$\Tr_{q^2/q}(\omega^{a+\frac{q+1}{2}\epsilon})\in C_0^{(2,q)}$, i.e., 
$a\in I_2-\frac{q+1}{2}\epsilon$. On the other hand, 
taking trace of both sides of  $\omega^{a+\frac{q+1}{2}\delta-c}=s\omega^{\frac{q+1}{2}\delta-c}+t\omega^{\frac{q+1}{2}\delta}$, we have 
$\Tr_{q^2/q}(\omega^{a+\frac{q+1}{2}\delta-c})=
s\Tr_{q^2/q}(\omega^{\frac{q+1}{2}\delta-c})+t\Tr_{q^2/q}(\omega^{\frac{q+1}{2}\delta})$. Since $\Tr_{q^2/q}(\omega^{\frac{q+1}{2}\delta})=0$ and $\Tr_{q^2/q}(\omega^{\frac{q+1}{2}\delta-c})\in C_0^{(2,q)}$, we obtain $\Tr_{q^2/q}(\omega^{a+\frac{q+1}{2}\delta-c})\in C_0^{(2,q)}$, i.e.,  $a\in I_2+c-\frac{q+1}{2}\delta$. Thus, $a\in (I_2-\frac{q+1}{2}\epsilon)\cap  (I_2+c-\frac{q+1}{2}\delta)$, and hence $G_3\subseteq Y_{3,c}$. Noting that $|G_3|=|Y_{3,c}|$, it follows that $G_3=Y_{3,c}$. 
Furthermore, since  $G_2\equiv G_3+(q+1)\,(\mod{2(q+1)})$ and 
$G_5=\{0,1,\ldots,2q+1\}\setminus (G_3\cup G_4 \cup \{0,c,q+1,c+q+1\})$, 
we have $G_4=Y_{4,c}$ and  $G_5=Y_{5,c}$

Finally, the character values of $R_4$ and $R_5$ are determined as  
$\psi_{\F_{q^2}}(\omega^a R_4)=\psi_{\F_{q^2}}(\omega^{a+q+1} R_3)$ and 
$\psi_{\F_{q^2}}(\omega^a R_5)=-1-\sum_{i=1}^4\psi_{\F_{q^2}}(\omega^a R_i)$. This completes the proof of the proposition. 
\qed

\begin{remark}
By the definition of $Y_{i,c}$, $i=1,2$, in Proposition~\ref{prop:asso}, it is clear that 
$-Y_{i,c}+c\equiv Y_{i,c}\,(\mod{2(q+1)})$, that is, the property~(P8). 
\end{remark}
Next, we show that the $Y_{i,c}$'s have property~(P9). 
\begin{proposition}
We have
\[
-(Y_{3,c}\cup Y_{4,c})+c\equiv Y_{5,c}\,(\mod{2(q+1)}). 
\]
\end{proposition}
\proof 
Since $Y_{i,c}=G_i$ for $i=3,4,5$ as in the proof of Proposition~\ref{prop:asso}, 
we have 
\begin{align*}
Y_{3,c}\cup Y_{4,c}=&\, \{a\,(\mod{2(q+1)})\,|\,\omega^a=s+t\omega^c,(s,t)\in S\times S\mbox{ or }N\times N\}, \\
Y_{5,c}=&\,\{a\,(\mod{2(q+1)})\,|\,\omega^a=s+t\omega^c,(s,t)\in S\times N\mbox{ or }S\times N\}. 
\end{align*}
Assume that $a \in -(Y_{3,c}\cup Y_{4,c})+c$. 
There are some $s',t'\in \F_q$ such that $\omega^a=s'+t'\omega^c$. 
On the other hand, since $a \in -(Y_{3,c}\cup Y_{4,c})+c$,  
$\omega^{-a+c}=s+t\omega^c$ for some $s,t\in S\times S$ or $N\times N$. Then, we have  
\begin{equation}\label{eq:stst}
(s\omega^{-c}+t)(s'+t'\omega^c)=1. 
\end{equation}
By multiplying both sides of \eqref{eq:stst} by $\omega^{\frac{q+1}{2}\epsilon}$ and taking trace, 
we have
\begin{equation}\label{eq:stst2}
ss'\Tr(\omega^{-c+\frac{q+1}{2}\epsilon})+(ts'+st')\Tr_{q^2/q}(\omega^{\frac{q+1}{2}\epsilon})+tt'\Tr_{q^2/q}(\omega^{c+\frac{q+1}{2}\epsilon})=\Tr_{q^2/q}(\omega^{\frac{q+1}{2}\epsilon}). 
\end{equation}
Since $\Tr_{q^2/q}(\omega^{\frac{q+1}{2}\epsilon})=0$ by the definition of $X_{1,c}$, \eqref{eq:stst2} is reduced to
\[
-ss' u\omega^{\frac{q+1}{2}(\epsilon-\delta)}=tt'v, 
\]
where $u=\Tr_{q^2/q}(\omega^{c+\frac{q+1}{2}\epsilon})$ and 
$v=\Tr_{q^2/q}(\omega^{-c+\frac{q+1}{2}\delta})$. Here, $u,v\in C_0^{(2,q)}$  by (\ref{eq:traceapp1}). Furthermore, $st^{-1}\in C_0^{(2,q)}$ by the definitions of $s,t$, and $-\omega^{\frac{q+1}{2}(\epsilon-\delta)}\in C_1^{(2,q)}$ by the definitions of $\epsilon,\delta$. Hence, either 
$(s',t')\in C_0^{(2,q)}\times C_1^{(2,q)}$ or $C_1^{(2,q)}\times C_0^{(2,q)}$ holds 
by noting that $(s',t')=(0,0)$ is impossible. 
Therefore, $a\in Y_{5,c}$, i.e., $-(Y_{3,c}\cup Y_{4,c})+c\subseteq Y_{5,c}$, follows. 
Finally, since $|-(Y_{3,c}\cup Y_{4,c})+c|=|Y_{3,c}\cup Y_{4,c}|=|Y_{5,c}|$, 
the statement of the proposition follows. 
\qed
\vspace{0.3cm}

Finally, we show that the $R_{i}'$'s have property~(P10). 
\begin{proposition}\label{re:dual}
Let $R_{i}'$, $i=1,2,3,4,5$, be defined as  in Subsection~\ref{subsec:set}.   
Then, $R_i'$, $i=1,2,3,4,5$, take the character values listed in Table~\ref{tab_2}. 
\end{proposition}
\proof
Since $Y_{i,c}-c+\frac{q+1}{2}\delta\equiv X_{i,c-\frac{q+1}{2}\delta+\frac{q+1}{2}\epsilon}$ by Lemma~\ref{lem:x1c}, Remark~\ref{rem:fives}~(5) and the definitions of $X_{i,c},Y_{i,c}$, $i=1,2,3,4,5$,  the statement follows from Proposition~\ref{prop:asso}. 
\qed

\end{document}